\documentclass[12pt]{article}
\usepackage[dvipdfmx]{graphicx}

\usepackage{amsmath,amssymb,amsthm, amscd}
\usepackage{ytableau,mathrsfs}
  
  \makeatletter
  \@addtoreset{equation}{section}
  \makeatother

\newcommand{\Z}{\mathbb{Z}}
\newcommand{\R}{\mathbb{R}}
\newcommand{\C}{\mathbb{C}}
\newcommand{\Q}{\mathbb{Q}}
\newcommand{\vphi}{\varphi}
\newcommand{\mcl}{\mathcal}
\newcommand{\mfk}{\mathfrak}
\newcommand{\mrm}{\mathrm}

\renewcommand{\tilde}[1]{\widetilde{#1}}
\newtheorem{dfn}{Definition}[section]
\newtheorem{thm}[dfn]{Theorem}
\newtheorem{lem}[dfn]{Lemma}
\newtheorem{prop}[dfn]{Proposition}
\newtheorem{cor}[dfn]{Corollary}
\theoremstyle{remark}
\newtheorem{rem}{Remark}[section]
\newtheorem{ex}{Example}[section]
\newtheorem{que}{Question}[section]

\begin{document}

\begin{center} 
{\bf {\LARGE Moduli of $G$-constellations and crepant resolutions I: the abelian case}}
\end{center} 
\vspace{0.4cm}

\begin{center}
{\large Ryo Yamagishi}
\end{center} 
\vspace{0.4cm}

\begin{abstract}
For a finite abelian subgroup $G\subset SL_n(\C)$, we study whether a given crepant resolution $X$ of the quotient variety $\C^n/G$ is obtained as a moduli space of $G$-constellations. In particular we show that, if $X$ admits a natural $G$-constellation family in the sense of Logvinenko over it with all fibers being indecomposable as $\C[\C^n]$-modules, then $X$ is isomorphic to the normalization of a fine moduli space of $G$-constellations.
\end{abstract} 

\renewcommand{\thefootnote}{\fnsymbol{footnote}} 

\section{Introduction}\label{1}

Let $G\subset SL_n(\C)$ be a finite subgroup. The aim of the present paper is to study how to realize crepant resolutions of $\C^n/G$ as moduli spaces of $G$-equivariant objects on $\C^n$. From the viewpoint of the McKay correspondence (cf. \cite{R}), it is expected that the geometry of a crepant resolution reflects the representation-theoretic properties of $G$, and giving a moduli description of a crepant resolution is very useful in studying such a phenomenon. Indeed, a derived category version of the McKay correspondence (called the {\it derived McKay correspondence}) is established in small dimensions by realizing crepant resolutions as certain moduli spaces (cf. \cite{BKR}).

A typical moduli space, as a candidate for a crepant resolution of $\C^n/G$, is the $G$-{\it Hilbert scheme} $G\text{-Hilb}$. It parametrizes $G$-{\it clusters} of $\C^n$, that is, $G$-invariant closed subschemes $Z\subset \C^n$ such that $H^0(\mcl{O}_Z)$ is isomorphic to the regular representation $R$ of $G$ as a $G$-module. It is known that $G\text{-Hilb}$ gives a projective crepant resolution of the quotient variety $\C^n/G$ when $n=2$ or $3$ \cite{IN}, \cite{BKR}. In higher dimensions, however, $G\text{-Hilb}$ usually behaves badly: it may not be irreducible nor normal in general. Even if we take the normalization of the irreducible component of $G\text{-Hilb}$ which is birational to $\C^n/G$, it can be singular or discrepant.

Nevertheless, we still have candidates of moduli spaces for crepant resolutions of $\C^n/G$, namely, moduli spaces of $G$-{\it constellations}. These moduli spaces were introduced as moduli of quiver representations by Kronheimer when $n=2$ \cite{Kr} and generalized to higher dimensions (cf. \cite{CI}). By definition, a $G$-constellation is a $G$-equivariant coherent sheaf $\mcl{F}$ on $\C^n$ such that $H^0(\mcl{F})$ is isomorphic to $R$ as a $G$-module. The structure sheaf of a $G$-cluster of $\C^n$ is by definition an example of a $G$-constellation. In order to consider a moduli space of $G$-constellations, we usually choose a {\it stability condition} $\theta$ which lives in a vector space
$$\Theta:=\{\theta\in\mrm{Hom}_\mathbb{Z}(R(G),\R)\mid \theta(R)=0\}$$
where $R(G)$ is the representation ring of $G$. Then, by using the notion of $\theta$-(semi)stable $G$-constellations introduced by King \cite{Kin}, we can construct a coarse moduli space $\mcl{M}_\theta$ of $\theta$-semistable $G$-constellations (up to S-equivalence) using the geometric invariant theory (GIT). We denote by $\bar{\mfk{M}}_\theta$ the normalization of the irreducible component (called the {\it coherent component}) of $\mcl{M}_\theta$ which is birational to $\C^n/G$. We obtain different moduli spaces by varying stability conditions $\theta$, and accordingly the vector space $\Theta$ admits a wall-and-chamber structure consisting of so-called GIT-chambers.

Besides the moduli spaces $\mcl{M}_\theta$ (or $\bar{\mfk{M}}_\theta$), it is also possible to consider moduli spaces which are not constructed by using the GIT construction above. For example, in \cite{L2} Logvinenko considered a family $\mcl{F}$ of $G$-constellations on a non-projective crepant resolution $X\to \C^3/G$ for certain $G$ and showed that $X$ can be regarded as a fine moduli space of $G$-constellations having $\mcl{F}$ as its universal family. Such a moduli space cannot be obtained as $\mcl{M}_\theta$ (or $\bar{\mfk{M}}_\theta$) since GIT quotients are always quasi-projective. 

The above construction of the non-projective moduli space is done by using the notion of a {\it geometrically natural family} ({\it gnat family} for short) introduced in \cite{L1},\cite{L3}. Every crepant resolution $X$ is birational to $\C^n/G$ by definition and hence admits a dense open subset parametrizing $G$-constellations. A gnat family $\mcl{F}$ on $X$ is defined as a family of $G$-constellations over $X$ which extends the family on this open subset (see Subsection \ref{4.1} for the precise definition). 

Giving a gnat family over $X$ is always possible but this is in fact not sufficient for $X$ to be a (fine) moduli space. In this paper we will construct a variety $\bar{\mfk{M}}_{X,\mcl{F}}$ for each gnat family $\mcl{F}$ which is a natural candidate for a normalized moduli space of $G$-constellations parametrized by $\mcl{F}$. $\bar{\mfk{M}}_{X,\mcl{F}}$ will be defined as a certain toric variety whose associated fan is determined by $X$ and $\mcl{F}$ (Definition \ref{candidate}). We will also construct a birational morphism $\alpha:X\to \bar{\mfk{M}}_{X,\mcl{F}}$ so that $\alpha$ is an isomorphism when $X$ is a fine moduli space of $G$-constellations with $\mcl{F}$ its universal family.

The main result of this article is to give a necessary and sufficient condition for $X$ to be such a moduli space.

\begin{thm}\label{main}(=Theorem \ref{main'})
Let $G\subset SL_n(\C)$ be a finite abelian subgroup, $X\to \C^n/G$ a (possibly non-projective) crepant resolution, and $\mcl{F}$ a gnat family over $X$. Then the morphism $\alpha:X\to \bar{\mfk{M}}_{X,\mcl{F}}$ is an isomorphism to the normalization of a fine moduli space of $G$-constellations whose universal family is pulled back by $\alpha$ to $\mcl{F}$ if and only if every $G$-constellation in $\mcl{F}$ is indecomposable (i.e., there is no nontrivial decomposition as a $\C[\C^n]$-module).
\end{thm}

We will see in Subsection \ref{4.3} that the fibers of $\mcl{F}$ are always indecomposable outside a codimension two subset $Z\subset X$ but a fiber $F_z$ for $z\in Z$ can be decomposable in general (cf. Example \ref{A3}). If the equivalent conditions in the theorem above are satisfied, we expect that the universal family $\mcl{F}$ induces an equivalence of the derived categories $D^b(X)$ and $D_G^b(\C^n)$ (cf. Question \ref{question}).

In this paper we will also show that a certain type of crepant resolutions (or more generally, relative minimal models) $X$ admits a description as a moduli space of $G$-constellations. More precisely, if $X$ is obtained by blowing-up $\C^n/G$ successively, then $X$ is isomorphic to the normalization of a coarse moduli space of $G$-constellations (Proposition \ref{prop:special}). The idea of the proof is almost due to the arguments in \cite{DLR} and \cite{J}, but we will rephrase these arguments in terms of gnat families.

This paper is organized as follows. In Section \ref{2} we review the construction of the moduli space $\bar{\mfk{M}}_\theta$. In Section \ref{3} we give a toric description of a crepant resolution (or a relative minimal model) $X$ of $\C^n/G$ and also review the quotient construction of toric varieties due to Cox as a preparation for the subsequent sections. In Section \ref{4} we explain the notion of a natural $G$-constellation family and introduce its related notions. We will also construct the morphism $\alpha:X\to \bar{\mfk{M}}_{X,\mcl{F}}$ in this section. In Section \ref{5}, we will prove Theorem \ref{main}. In Section \ref{6}, we show how to realize $X$ as a coarse moduli space of $G$-constellations when $X$ is obtained by blowing-ups.

\vspace{5mm}
\noindent{\underline{Conventions \& Notations}\\
In this paper a {\it variety} means an integral separated scheme of finite type over $\C$. By a {\it point} of a variety, we mean a closed point. We always assume that an action of an algebraic group on a variety is algebraic. For a group $H$, we denote by $\chi(H)$ the character group of $H$. For an abelian group $A$, we denote by $A_\R$ its scalar extension $A\otimes_\Z \R$. For a finite dimensional $\R$-vector space $W$ and a subset $S\subset W$, we denote by $\mrm{Cone}(S)$ the cone generated by $S$.

\vspace{5mm}
\noindent{\bf Acknowledgements}\\
This work was supported by World Premier International Research Center Initiative (WPI), MEXT, Japan, and by JSPS KAKENHI Grant Number JP19K14504.

\section{Moduli spaces of $G$-constellations for abelian $G$}\label{2}
In this section we recall the notion of a $G$-constellation for a finite abelian subgroup $G\subset SL_n(\C)$ and review the construction of a moduli space of $G$-constellations using the McKay quiver associated to $G$. See \cite{CMT} for a more detailed exposition and concrete examples.

The action of $G$ on $V:=\C^n$ naturally induces an action on the coordinate ring $\C[V]$ so that
\begin{equation}\label{action}
(g\cdot f)(v)=g^{-1}\cdot (f(v))
\end{equation}
for any $g\in G, f\in\C[V]$ and $v\in V$. If a $\C[V]$-module $M$ admitting a linear $G$-action satisfies
$$g\cdot (fm)=(g\cdot f)m$$ 
for any $g\in G, f\in\C[V]$ and $m\in M$, then we say that $M$ is {\it $G$-equivariant}.

\begin{dfn}(cf. \cite{CI})
A $G$-equivariant $\C[V]$-module $M$ is called a {\bf $G$-constellation} if it is isomorphic to the regular representation $R$ of $G$ as a $G$-module.
\end{dfn}

\begin{rem}
Usually a $G$-constellation on a $G$-scheme $S$ is defined as a $G$-equivariant coherent sheaf $F$ on $S$ such that $H^0(F)\cong R$ as a $G$-module (cf. \cite{CI}). In this paper we only consider $G$-constellations on the affine space $V$ and thus we regard them as modules rather than sheaves.
\qed
\end{rem}

Let $\mrm{Irr}(G)=\{\rho_0,\dots,\rho_{r-1}\}\,(r=|G|)$ be the set of isomorphism classes of the irreducible representations of $G$ with $\rho_0$ being the trivial representation. Since every $\rho_i$ is 1-dimensional, we will often regard each $\rho_i$ as a character $G\to \C^*$ of $G$ and vice versa without mentioning it. Note that we have $R\cong\bigoplus_{i=0}^r \rho_i$. We choose coordinates $x_1,\dots,x_n$ of $V$ so that every element $g\in G$ acts on each $x_j$ as a scalar, that is, there is a character $\chi_j: G\to \C^*$ of $G$ such that $g\cdot x_j=\chi_j(g)x_j$. We fix such $x_j$'s hereafter.

Next we introduce the McKay quiver $Q$ associated to $G$. The quiver $Q$ is the oriented graph consisting of the vertices $v_i=v_{\rho_i}\,(i=0,\dots,r-1)$ indexed by the representations $\rho_i$ and of arrows $\alpha_{i,j}\,(i=0,\dots,r-1,\,j=1,\dots,n)$ each of which has its tail and head at $v_i$ and $v_{\rho_i\otimes \chi_j}$ respectively. We denote by $k(i,j)$ the integer such that $v_{k(i,j)}=v_{\rho_i\otimes \chi_j}$.

We construct a moduli space parametrizing isomorphism classes of $G$-constellations using the quiver $Q$. Let $A=\C[x_{i,j}]_{0\le i\le r-1,1\le j\le n}$ be the polynomial ring with $rn$ variables $x_{i,j}$. For each triple $(i,j,j')$ of integers $i\in\{0,\dots,r-1\}$ and $j,j'\in\{1,\dots, n\}$ with $j\ne j'$, we define $f_{i,j,j'}\in A$ as the binomial 
$$x_{k(i,j),j'}x_{i,j}-x_{k(i,j'),j}x_{i,j'}.$$
Let $I\subset A$ be the ideal generated by all $f_{i,j,j'}$'s, and set $\mcl{R}=\mrm{Spec}\,A/I$.

We can regard $\mcl{R}$ as a space of $G$-constellations as follows. Let us consider the vector space $\C^r$ with the standard basis $e_0,\dots,e_{r-1}$. We regard $\C^r$ as the regular $G$-representation by letting $G$ act on $\C e_i$ as the 1-dimensional representation $\rho_i$. For any closed point $p=(p_{i,j})\in \mcl{R}\subset \C^{rn}$, we endow $\C^r$ with a structure of a $\C[V]$-module by setting
\begin{equation}\label{module}
x_j\cdot e_i=p_{i,j}e_{k(i,j)}.
\end{equation}
The relations $p_{k(i,j),j'}p_{i,j}-p_{k(i,j'),j}p_{i,j'}=0$ imply the commutativity of the actions of $x_j$'s, and thus this $\C[V]$-module structure is well-defined. One can also check easily that this module is $G$-equivariant.

Conversely, if we once fix a $G$-module isomorphism between the regular representation $R=\bigoplus_{i=0}^{r-1}\C e_i$ and a given $G$-constellation $M$, we can associate a closed point $p=(p_{i,j})$ of $\mcl{R}$ by the equation (\ref{module}). Note that the $G$-equivariance of $M$ guarantees that $x_j(\C e_i)\subset\C e_{k(i,j)}$. Since $\C^r$ is a direct sum of 1-dimensional representations, its $G$-module automorphism group is the algebraic torus $T^r:=(\C^*)^{\times r}$ whose $i$-th component acts as scaler multiplication on $\C e_i$. Via the correspondence between $\C[V]$-module structures on $\C^r$ and points of $\mcl{R}$, we obtain a $T^r$-action on $\mcl{R}$ given by
\begin{equation}\label{Tr}
(t_0,\dots,t_{r-1})\cdot (p_{i,j})_{i,j}=(t_i^{-1}t_{k(i,j)}p_{i,j})_{i,j}
\end{equation}
for any $(t_1,\dots,t_r)\in T^r$ and $(p_{i,j})_{i,j}\in \mcl{R}$. Then we can conclude that the set of isomorphism classes of $G$-constellations are in one-to-one correspondence with the set of $T^r$-orbits of $\mcl{R}$. Note that the diagonal subgroup $\Delta:=\{(t,t,\dots,t)\}\subset T^r$ acts trivially on $\mcl{R}$. We put $\bar{T}^{r-1}=T^r/\Delta$.

Since the orbit space $\mcl{R}/T^r(=\mcl{R}/\bar{T}^{r-1})$ does not have a reasonable structure as a scheme, we usually introduce the notion of a stability condition coming from the geometric invariant theory in order to construct a moduli space of $G$-constellations. In our case a character $\chi$ of $\bar{T}^{r-1}$ gives a $\bar{T}^{r-1}$-action on the total space of the trivial line bundle $\mcl{R}\times\C$ by
$$\bar{T}^{r-1}\times(\mcl{R}\times\C)\to\mcl{R}\times\C;(\mathbf{t},(p,a))\mapsto (\mathbf{t}\cdot p,\chi(\mathbf{t})^{-1}a).$$
Note that this action is the same as the natural one coming from the $\bar{T}^{r-1}$-action on $\mcl{R}$ (as in (\ref{action})) if $\chi$ is trivial. Then the set of $\bar{T}^{r-1}$-invariant global sections is the same as the $\C$-vector space
$$(R/I)_\chi=\{f\in R/I\mid \mathbf{t}\cdot f=\chi(\mathbf{t})f\}.$$
The notion of semistable and stable points of $\mcl{R}$ with respect to $\chi$ are defined as follows (cf. \cite{MFK}).

\begin{dfn}
Let $\chi:\bar{T}^{r-1}\to\C^*$ be a character. We say that a point $p\in\mcl{R}$ is $\chi${\bf -semistable} if $f(p)\ne0$ for some $f\in (R/I)_{\chi^k}$ with $k\ge1$. If furthermore the $\bar{T}^{r-1}$-orbit $\bar{T}^{r-1}\cdot p$ is closed in the open subset $\{f\ne0\}\subset \mcl{R}$ and is of dimension $\dim \bar{T}^{r-1}=r-1$, we say that $p$ is $\chi${\bf -stable}. We denote the set of $\chi$-semistable (resp. $\chi$-stable) points of $\mcl{R}$ by $\mcl{R}^{\chi\text{-ss}}$ (resp. $\mcl{R}^{\chi\text{-s}}$). If we have $\mcl{R}^{\chi\text{-ss}}=\mcl{R}^{\chi\text{-s}}$, we say that $\chi$ is {\bf generic}.
\end{dfn}

Note that $\chi$ and its positive power give the same (semi)stable locus. In our situation where $\mcl{R}$ is a space of quiver representations, it was observed by King \cite{Kin} that we can interpret the notion of $\chi$-(semi)stability in a representation-theoretic way as follows. For a given character $\chi:\bar{T}^{r-1}\to\C^*$, we can write $\chi(\overline{t_0,\dots,t_r})=t_0^{\theta_0}\dots t_{r-1}^{\theta_{r-1}}$ for some $\theta_i\in\mathbb{Z}$. Then the assignment
$$\begin{aligned}
\chi(&\bar{T}^{r-1}) & \to& & \bar{\Theta}:=\{&\theta\in\mrm{Hom}_\mathbb{Z}(R(G),\Z)\mid \theta(R)=0\}\\
&\chi & \mapsto&  & &\theta_\chi=(\rho_i\mapsto \theta_i,\;\forall i)
\end{aligned}$$
is an isomorphism of abelian groups. King showed that a point $p\in\mcl{R}$ is $\chi$-semistable (resp. $\chi$-stable) if and only if the corresponding $\C[V]$-module to $p$ is $\theta_\chi$-semistable (resp. $\theta_\chi$-stable) \cite[Proposition 3.1]{Kin}. Here the (semi)stability of a module is defined as follows.

\begin{dfn}(cf. \cite{Kin})
Let $M$ be a $\C[V]$-module. For an element $\theta\in\bar{\Theta}$, we say that $M$ is $\theta${\bf -semistable} if $\theta(M)=0$ and $\theta(M')\ge0$ for all submodule $M'\subset M$. Such an $M$ is called $\theta${\bf -stable} if the only submodules $M'\subset M$ satisfying $\theta(M')=0$ are $0$ and $M$.
\end{dfn}

We will use the scaler extension 
$$\Theta:=\{\theta\in\mrm{Hom}_\mathbb{Z}(R(G),\R)\mid \theta(R)=0\}=\bar{\Theta}\otimes_\Z \R$$
as the space of stability conditions rather than $\bar{\Theta}$ itself in order to consider polyhedral cones inside $\Theta$ (or $\chi(\bar{T}^{r-1})_\R$). In fact $\theta$-stability makes sense also for elements of $\Theta$ (see Remark \ref{rem:stability} below).

By considering the (semi)stable locus of $\mcl{R}$ instead of $\mcl{R}$ itself, we obtain the quotient space as a variety (after taking an irreducible component). To be more precise, let us recall the notion of good and geometric quotients.

\begin{dfn}
Let $H$ be a reductive group acting on a variety $Y$. An $H$-invariant morphism $p:Y\to Z$ to a variety is a {\bf good quotient} if it is affine and the induced map $\mcl{O}_Z\to (p_*\mcl{O}_Y)^H$ to the $H$-invariant part of $p_*\mcl{O}_Y$ is an isomorphism. If additionally every fiber of $p$ is a single $H$-orbit, we call $p$ a {\bf geometric quotient}.
\end{dfn}

\begin{rem}
Later we will also consider categorical quotients (in the category of algebraic varieties) in Subsection \ref{4.3}. By definition, an $H$-invariant morphism $p:Y\to Z$ to a variety $Z$ is a {\bf categorical quotient} if for any $H$-invariant morphism $p':Y\to Z'$ to a variety, there exists a morphism $f:Z\to Z'$ such that $p'=f\circ p$. It is known that a good quotient is automatically a categorical quotient (cf. \cite{MFK}).
\qed
\end{rem}

\vspace{3mm}

By geometric invariant theory (cf. \cite{MFK}), we can form a good quotient
$$\pi_\chi:\mcl{R}^{\chi\text{-ss}}\to \mcl{R}/\!/_\chi \bar{T}^{r-1}:=\mrm{Proj}\left(\bigoplus_{k=0}^\infty (R/I)_{\chi^k}\right)$$
whose restriction $\mcl{R}^{\chi\text{-s}}\to \pi_\chi(\mcl{R}^{\chi\text{-s}})$ is a geometric quotient. For $\theta\in\Theta$, now we can define the moduli space of $\theta$-semistable $G$-constellations $\mcl{M}_\theta$ as $\mcl{R}^{\chi\text{-ss}}/\!/_\chi \bar{T}^{r-1}$ with the character $\chi\in\chi(\bar{T}^{r-1})_\R$ corresponding to $\theta$ via $\chi(\bar{T}^{r-1})_\R\cong \Theta$. It is known that $\mcl{M}_\theta$ is a coarse moduli space parametrizing S-equivalence classes of $\theta$-semistable $G$-constellations \cite[Proposition 5.2]{Kin} where two representations of a quiver are {\it S-equivalent} if they have the same composition factors. If $\theta$ is generic, then it follows from \cite[Proposition 5.3]{Kin} that $\mcl{M}_\theta$ is a fine moduli space parametrizing isomorphism classes of $\theta$-stable $G$-constellations.

\vspace{3mm}

\begin{rem}\label{rem:stability}
The $\chi$-(semi)stability can make sense even if $\chi$ is not on the lattice $\chi(\bar{T}^{r-1})$ since the GIT-equivalence classes gives a wall-and-chamber structure on $\chi(\bar{T}^{r-1})$ where two characters of $\bar{T}^{r-1}$ are GIT-equivalent if they give the same semistable locus (cf. \cite[\S2]{T}). The space $\Theta$ also admits a wall-and-chamber structure through the identification $\chi(\bar{T}^{r-1})_\R\cong\Theta$, and we see from the definition of $\theta$-(semi)stability that each wall of $\Theta$ is realized as a part of a hyperplane $H_S=\{\sum_{i\in S}\theta(\rho_i)=0\}\subset \Theta$ for a proper subset $S\subset \{0,\dots,r-1\}$. Note that the wall-and-chamber structure is not necessarily given by a hyperplane arrangement.
\qed
\end{rem}

\vspace{3mm}

Let $\mathbf{0}\in\chi(\bar{T}^{r-1})$ be the trivial character. Then for any character $\chi\in\chi(\bar{T}^{r-1})$, we have a natural morphism from $\mcl{R}/\!/_\chi \bar{T}^{r-1}\to \mcl{R}/\!/_\mathbf{0} \bar{T}^{r-1}$ which factors through a closed immersion $\C^n/G\hookrightarrow \mcl{R}/\!/_\mathbf{0} \bar{T}^{r-1}$ \cite[Proposition 2.2]{CI}. The affine variety $\mcl{R}$ is not irreducible in general \cite[Example 3.1]{CMT}, but it always has a unique ($\bar{T}^{r-1}$-invariant) irreducible component $\mcl{V}\subset\mcl{R}$ which does not lie on any
coordinate hyperplane \cite[Theorem 3.10]{CMT}. The component $\mcl{V}$ also satisfies $\mcl{V}/\!/_\mathbf{0} \bar{T}^{r-1}= \C^n/G$ \cite[Proposition 4.1]{CMT}. Accordingly the moduli space $\mcl{M}_\theta$ also has a distinguished irreducible component $\mfk{M}_\theta$, which is called the {\it coherent component}. $\mfk{M}_\theta$ is nonnormal in general but is always birational to $\C^n/G$ via the natural map. Thus, the normalization $\bar{\mfk{M}}_\theta$ of $\mfk{M}_\theta$ is a candidate of a projective crepant resolution of $\C^n/G$ (if it exists).

\vspace{3mm}

\begin{rem}\label{nonproj1}
The moduli space $\bar{\mfk{M}}_\theta$ for generic $\theta$ is always projective over $\C^n/G$ \cite[Proposition 2.2]{CI}. However, we can also construct a non-projective (fine) moduli space by choosing a suitable $\bar{T}^{r-1}$-invariant open subset of $\mcl{R}$ not coming from a character of $\bar{T}^{r-1}$ (cf. Remark \ref{nonproj2}).
\qed
\end{rem}

\section{Crepant resolutions of abelian quotient singularities}\label{3}

In this section we review the basic facts about abelian quotient singularities and their crepant resolutions. In this paper we describe these objects as toric varieties. So we first recall basic notions about toric varieties. See e.g. \cite{Co2} for a comprehensive exposition for toric varieties.

\subsection{Basics of toric varieties}\label{3.1}
Let $N=\mathbb{Z}^n$ be the free abelian group of rank $n$. In this paper, a subset $\sigma\subset N_\mathbb{R}=N\otimes_\mathbb{Z}\mathbb{R}$ is called a {\it cone} if there are finitely many elements $v_1,\dots,v_s\in N$ such that $\sigma=\sum_{i=1}^s\mathbb{R}_{\ge0} v_i$ and if $\sigma\cap (-\sigma)=\{0\}$. Such $v_i$'s are called {\it generators} of $\sigma$. A cone $\sigma$ is {\it smooth} if we can take generators of $\sigma$ as a part of basis of $N$. As a weaker notion, we say that a cone $\sigma$ is {\it simplicial} if the number of minimal generators of it is the same as the dimension of $\sigma$. Here, the dimension of a cone $\sigma$ is defined as the smallest integer $k$ such that $\sigma$ is contained in a $k$-dimensional linear subspace of $N_\mathbb{R}$. 

Let $M=\mrm{Hom}_\mathbb{Z}(N,\mathbb{Z})$ be the dual abelian group of $N$ and let $f_1,\dots,f_n\in M$ be the dual basis of the standard basis $e_1,\dots,e_n\in N$. We will often identify the group algebra $\C[M]=\bigoplus_{m\in M} \chi^m$ with the Laurent polynomial ring $\C[x_1^{\pm1},\dots,x_n^{\pm1}]$ by setting $x_i=\chi^{f_i}$. For a cone $\sigma\subset N_\mathbb{R}$, the dual cone of $\sigma$ is defined as
$$\sigma^\vee=\{m\in M_\mathbb{R} \mid \langle u,m \rangle\ge0,\,\forall u\in\sigma\}\subset M_\mathbb{R}$$
where $\langle -,- \rangle: N_\mathbb{R}\times M_\mathbb{R}\to \mathbb{R}$ is the natural pairing. The semigroup $\sigma^\vee\cap M$ is finitely generated by Gordan's lemma, and we define the affine toric variety $X(\sigma)$ for $\sigma$ as the spectrum $\mrm{Spec}\,\C[\sigma^\vee\cap M]$. As is well known, the variety $X(\sigma)$ is smooth if and only if the cone $\sigma$ is smooth (cf. \cite[1.3.12]{Co2}). Moreover, $X(\sigma)$ is $\Q$-factorial (i.e., some positive multiple of any Weil divisor is Cartier) if and only if the cone $\sigma$ is simplicial \cite[Proposition 4.2.7]{Co2}.

Next we consider a toric variety associated to a fan in $N_\mathbb{R}$. A {\it fan} $\Sigma$ in $N_\mathbb{R}$ is a collection $\{\sigma\}$ of cones in $N_\mathbb{R}$ satisfying the following conditions:
\begin{itemize}
\item every face of any cone $\sigma\in \Sigma$ again belongs to $\Sigma$
\item the intersection of any two cones $\sigma,\tau\in \Sigma$ again belongs to $\Sigma$
\end{itemize}
where a {\it face} of a cone $\sigma$ is the intersection of $\sigma$ and the linear subspace $H_m=\{u\in N_\mathbb{R} \mid \langle u,m \rangle=0\}$ of $N_\mathbb{R}$ for some $m\in \sigma^\vee$ such that $\sigma\subset H^+_m=\{u\in N_\mathbb{R} \mid \langle u,m \rangle\ge0\}$. Note that a face of a cone is again a cone. For a $k$-dimensional cone $\sigma$, a $(k-1)$-dimensional face of $\sigma$ is called a {\it facet} of $\sigma$. A 1-dimensional cone is called a {\it ray}. For a fan $\Sigma$, we also define $\Sigma_{\mrm{max}}$ as the set of cones in $\Sigma$ which are maximal with respect to inclusion.

For a fan $\Sigma$ in $N_\mathbb{R}$, we can associate to it a variety $X(\Sigma)$ which has an open covering by $X(\sigma)$'s for $\sigma\in\Sigma$. This variety $X(\Sigma)$ is called the toric variety associated to the fan $\Sigma$. We do not give the precise construction of $X(\Sigma)$ because we will not use it later. The important feature of $X(\Sigma)$ is that it admits an action by the algebraic torus $T_N=N\otimes_\mathbb{Z}\C^*$ and that the open subset $X(\sigma)\subset X(\Sigma)$ for each $\sigma\in\Sigma$ is $T_N$-invariant and the induced $T_N$-action on $X(\sigma)$ coincides with the natural one, namely, the one coming from the $T_N$-action
$$T_N\times \C[M]\to\C[M]; (u\otimes t,\chi^m)\mapsto t^{\langle u,m \rangle}\chi^m$$
on $\C[M]$. In this paper we sometimes call $T_N$ the {\it big torus} of of the toric variety $X(\Sigma)$.

\subsection{Abelian quotient singularities and toric partial resolutions}\label{3.2}
In this subsection we briefly review the construction of abelian quotient singularities and their partial resolutions as toric varieties. For further details, we refer the reader to \cite{DHZ}. 

Let $G\subset SL_n(\C)$ be a finite abelian subgroup. We choose coordinates $x_1,\dots,x_n$ of $V=\C^n$ as in section \ref{2}. The quotient variety $\C^n/G$ is defined as the spectrum $\mrm{Spec}\,\C[V]^G$ of the $G$-invariant subring of $\C[V]=\C[x_1,\dots,x_n]$. We first realize $\C^n/G$ as an affine toric variety.

By the choice of $x_j$'s, we can uniquely write
\begin{equation}\label{age}
g\cdot x_j=e^{\frac{2\pi ia_j}{r_g}} x_j,\;0\le a_j<r_g
\end{equation}
for any $j$ and $g\in G$ with its order $r_g$. Using this, we associate an element $v_g=1/r_g(a_1,\dots,a_n)\in \Q^n$ to each $g\in G$. Regarding $N=\mathbb{Z}^n$ as a subset of $\Q^n$, we define an abelian subgroup $N_G$ of $\Q^n$ as the subgroup generated by $N=\mathbb{Z}^n$ and $\{v_g \mid g\in G\}$. Then we have the following exact sequence of abelian groups
\begin{equation}\label{seq1}
0\to N\stackrel{\iota}{\to} N_G\to G\to 0
\end{equation}
where the homomorphism $N_G\to G$ sends $v=\sum_{g\in G}n_g v_g+v'\,(v'\in N)$ to $\prod_{g\in G} g^{n_g}$. Note that this is independent of the choice of the expression of $v$, and then the exactness of the sequence (\ref{seq1}) is clear. As for the dual side, let $M_G=\mrm{Hom}_\mathbb{Z}(N_G,\mathbb{Z})$ be the dual abelian group of $N_G$. Then, similarly we have the following exact sequence of abelian groups
\begin{equation}\label{seq2}
0\to M_G\to M\to \chi(G)\to 0
\end{equation}
where $\chi(G)$ is the group of characters of $G$. The homomorphism $M\to \chi(G)$ is given by sending $m\in M$ to the character $\chi_0$ such that $g\cdot \chi^m=\chi_0(g)\chi^m$ where we use the identification $\chi^{f_j}=x_j$ as in the previous subsection. The surjectivity of this homomorphism follows from the effectiveness of the action of $G$ on $V$ (and hence on $\C[M]$). The exactness in the middle of the sequence (\ref{seq2}) also follows since $g\cdot \chi^m=\chi^m$ for $m=\sum_{j=1}^n b_j f_j$ if and only if the inner product of $v_g$ and $(b_1,\dots,b_n)$ is in $\mathbb{Z}$. Therefore, the algebra $\C[M_G]$ is exactly the $G$-invariant part of $\C[M]$.

The vector space $V=\C^n$ is identified with the toric variety $X(\sigma_0)$ for the cone $\sigma_0=\sum_{j=1}^n \mathbb{R}_{\ge0} e_j\subset N_\mathbb{R}$. The quotient $\C^n/G$ is also realized as the toric variety $X(\sigma_0^G)$ where the cone $\sigma_0^G$ is also generated by the basis $\{e_i\}$ of $N$ but is regarded as the cone inside $(N_G)_\mathbb{R}$. Note, in particular, that the associated big tori $T_N$ and $T_{N_G}$ for $X(\sigma_0)$ and $X(\sigma_0^G)$ respectively are the different (but isomorphic) ones satisfying $T_{N_G}/T_N\cong N_G/N\cong G$.

We introduce the notion of {\it age} following \cite{IR}.

\begin{dfn}\label{junior}
For each element $g\in G$, we consider the vector $v_g=1/r_g(a_1,\dots,a_n)\in \Q^n$ defined above. Then the {\bf age} of $g$ is the number defined as $\frac{1}{r_g}\sum_{j=1}^n a_i$, which is an integer since $g\in SL(V)$. We call an element $g\in G$ {\bf junior} if the age of $g$ is 1.
\end{dfn}

As we will see below, the junior elements of $G$ bijectively correspond to the exceptional prime divisors of a crepant resolution of $\C^n/G$.

Since $G$ is a subgroup of $SL_n(\C)$ and in particular contains no pseudo- (or complex) reflections, the result of Chevalley-Shephard-Todd (cf. \cite{Ch},\cite{ST}) shows that the quotient $\C^n/G$ is singular or, equivalently, the cone $\sigma_0^G\subset (N_G)_\mathbb{R}$ is not smooth unless $G$ is trivial. It is well-known that the singularity of a toric variety can be resolved in the category of toric varieties, where a morphism between two toric varieties $X(\Sigma_1)$ and $X(\Sigma_2)$ with fans $\Sigma_i\subset (N_i)_\mathbb{R}\,(i=1,2)$ in this category is given by a pair of a group homomorphism $\vphi:T_{N_1}\to T_{N_2}$ of the associated tori and a $\vphi$-equivariant morphism $X(\Sigma_1)\to X(\Sigma_2)$ of varieties. We call such a morphism a {\it toric morphism}. In terms of fans, giving a toric morphism $\vphi:X(\Sigma_1)\to X(\Sigma_2)$ is equivalent to giving a homomorphism $\bar{\vphi}:N_1\to N_2$ such that, for any cone $\sigma\in \Sigma_1$, there exists a cone $\tau\in \Sigma_2$ satisfying $\bar{\vphi}_\mathbb{R}(\sigma)\subset \tau$ where $\bar{\vphi}_\mathbb{R}:(N_1)_\mathbb{R}\to (N_2)_\mathbb{R}$ is the scalar extension of $\bar{\vphi}$. It is known that a toric morphism $\vphi:X(\Sigma_1)\to X(\Sigma_2)$ is proper if and only if the corresponding homomorphism $\bar{\vphi}:N_1\to N_2$ satisfies $\bar{\vphi}_\mathbb{R}^{-1}(|\Sigma_2|)=|\Sigma_1|$ where $|\Sigma_i|=\bigcup_{\sigma\in\Sigma_i}\sigma\subset (N_i)_\mathbb{R}$ is the support of $\Sigma_i$ \cite[Theorem 3.4.11]{Co2}.

Let $X(\sigma)$ be a (singular) affine toric variety. Using the facts above, we can conclude that giving a toric resolution of $X(\sigma)$ (i.e., a proper birational toric morphism to $X(\sigma)$ from a smooth toric variety) is equivalent to giving a subdivision of $\sigma$ into smooth cones (namely, giving a fan $\Sigma$ which consists of smooth cones and satisfies $|\Sigma|=\sigma$).

In this paper we are particularly interested in crepant resolutions of $\C^n/G=X(\sigma_0^G)$. By definition, a resolution $\pi:X\to \C^n/G$ is {\it crepant} if it preserves the canonical divisor i.e., $K_X=\pi^*K_{\C^n/G}$. Note that $\C^n/G$ is Gorenstein \cite{W} and thus this equation makes sense. If $\pi$ is a toric resolution obtained by a refinement $\Sigma$ of the cone $\sigma_0^G$, we have the following characterization of crepant resolutions in terms of the fan $\Sigma$ (see e.g. \cite[\S4]{DHZ}).

\begin{prop}\label{fan}
Using the notation as above, a toric resolution $\pi:X(\Sigma)\to \C^n/G$ is crepant if and only if every 1-dimensional cone $\sigma\in \Sigma$ is generated by an element $u=(u_1,\dots,u_n)\in N_G\subset \Q^n$ (i.e. $\sigma=\mathbb{R}_{\ge0}u$) such that $\sum_{j=1}^n u_j=1$.
\end{prop}

\begin{rem}\label{IR}
This proposition shows that the set of the exceptional prime divisors of $\pi$ is in one-to-one correspondence with the set of the junior elements in $G$ since each element $u=(u_1,\dots,u_n)\in N_G$ satisfying $\sum_{j=1}^n u_j=1$ and $0\le u_i<1$ is identified with a junior element via the map $N_G\to G$ in the sequence (\ref{seq1}). In fact a similar correspondence is known to hold for crepant resolutions of $\C^n/G$ with not-necessarily-abelian $G$. This was proved by Ito and Reid as a part of the (cohomological) McKay correspondence \cite{IR}.
\qed
\end{rem}

\vspace{3mm}

In dimension 2 or 3, a crepant resolution of $\C^n/G$ always exists. For $n\ge4$, crepant resolutions of $\C^n/G$ do not exist in general, or in other words there does not exist a subdivision of $\sigma_0^G$ into smooth cones satisfying the condition in Proposition \ref{fan}. However, if we allow simplicial cones instead of smooth cones, such a subdivision $\Sigma$ exists. A toric variety $X(\Sigma)$ obtained in this way is a so-called (relative) minimal model of $\C^n/G$. The precise definition is given as follows.

\begin{dfn}\label{minimal}
A proper birational morphism $\varphi:Y\to \C^n/G$ from a normal variety $Y$ is a {\bf relative minimal model} if it is crepant and $Y$ has at worst $\Q$-factorial terminal singularities. 
\end{dfn}

One can indeed show that a toric morphism $\varphi:X(\Sigma)\to \C^n/G$ obtained by a subdivision $\Sigma$ of $\sigma_0^G$ is a relative minimal model if and only if all cones of $\Sigma$ are simplicial and the set $\{v_g \mid g\in G:\text{junior}\}\cup\{e_1,\dots,e_n\}\subset N_G$ is exactly the set of the primitive generators of the rays of $\Sigma$ (cf. \cite[\S11.4]{Co2}). One can obtain such a fan $\Sigma$, for example, by applying star subdivisions to $\sigma_0^G$ (see Section \ref{6}). 

In the next subsection we will see that every relative minimal model (in particular every crepant resolution) of $\C^n/G$ is obtained as a toric variety.

\vspace{3mm}

\begin{rem}
In this paper we do not assume $\varphi$ is projective.
\qed
\end{rem}

\subsection{Quotient construction of toric varieties}\label{3.3}

In this subsection we briefly recall the construction of a toric variety as a quotient of an affine space, which is originally due to Cox \cite{Co1} (see also \cite[\S5.1]{Co2}). This will be used for the construction of the morphism $\alpha$ in Theorem \ref{main} (cf. Subsection \ref{4.3}).

Let $X=X(\Sigma)$ be a toric variety associated to a fan $\Sigma\subset N_\R$ with $N=\Z^n$, and let $\tau_1,\dots,\tau_m$ be the rays of $\Sigma$. If the dimension, say $d$, of the vector subspace $N'\subset N_\R$ generated by $\tau_k$'s is strictly smaller than $n$, then we also choose $n-d$ rays $\tau_{m+1},\dots,\tau_{m+n-d}$ of $N_\R$ so that $\tau_1,\dots,\tau_{m+n-d}$ generate $N_\R$. Let $S_X=\C[x_{\tau_1},\dots,x_{\tau_{m+n-d}}]$ be the polynomial ring whose variables are indexed by $\tau_k$'s. We define $T_\Sigma$ as the character group $\chi(\mrm{Cl}(X))$ of the divisor class group of $X$. Note that we have the following exact sequence
\begin{equation}\label{dclass}
0 \to M' \to \Z^m\to \mrm{Cl}(X) \to 0
\end{equation}
where $M'$ is the dual of $N'$. In particular, $\mrm{Cl}(X)$ is an abelian group of rank $m-d$ (cf. \cite[Lemma 5.1.1]{Co2}). Note also that the surjection $\Z^m\to \mrm{Cl}(X)$ in the sequence (\ref{dclass}) gives rise to an inclusion $T_\Sigma\subset (\C^*)^m$.

The affine space $\mrm{Spec}\,S_X$ is regarded as a toric variety with the lattice $\tilde{N}:=\Z^{m+n-d}$. We let $(\C^*)^m$ (and hence $T_\Sigma$) act on $\mrm{Spec}\,S_X$ by multiplication on the first $m$ coordinates of $S_X$. Let $\tilde{e}_1,\dots,\tilde{e}_{m+n-d}$ be the standard basis of $\tilde{N}$. We define a fan $\tilde{\Sigma}$ in $\tilde{N}_\R$ as
$$\tilde{\Sigma}=\{\mrm{Cone}(\tilde{e}_{i_1},\dots,\tilde{e}_{i_\ell}) \mid {\mrm{Cone}}(\tau_{i_1},\dots,\tau_{i_\ell})\in \Sigma\}.$$
Let $U_{\tilde{\Sigma}}\subset \mrm{Spec}\,S_X$ be the toric open subvariety defined by $\tilde{\Sigma}$. Then the surjection $\tilde{N} \to N$ which sends $\tilde{e}_i$ to the generator of $\tau_i$ induces a toric morphism $U_{\tilde{\Sigma}}\to X$.

\begin{prop}\label{quot}(cf. \cite[Theorem 5.1.11 and (5.1.11)]{Co2})
The morphism $U_{\tilde{\Sigma}}\to X$ is a good quotient of $U_{\tilde{\Sigma}}$ by the $T_\Sigma$-action. Moreover, $U_{\tilde{\Sigma}}\to X$ is a geometric quotient if and only if $X$ is $\Q$-factorial.
\end{prop}

\vspace{3mm}

Note that a crepant resolution (or a relative minimal model) $X$ of $\C^n/G$ is not a priori a toric variety and the polynomial ring $S_X$ is not defined. Nevertheless, we can still attach a graded ring $\mrm{Cox}(X)$ called the {\it Cox ring} of $X$. $\mrm{Cox}(X)$ is a finitely generated commutative ring graded by $\mrm{Cl}(X)$. We omit the precise definition as we will not use it. We refer the reader to \cite{ADHL} for an exposition. The only feature we will use is that $X$ is toric and $\mrm{Cox}(X)$ is isomorphic to $S_X$ as a graded ring if $\mrm{Cox}(X)$ is a polynomial ring. As we will see below, this is indeed the case.

For the given finite abelian subgroup $G\subset SL_n(\C)$, let $\{g_1,\dots,g_m\}$ be the set of junior elements of $G$. We define a $\mathbb{Z}^m$-grading on the polynomial ring $\C[V]=\C[x_1,\dots,x_n]$ as follows. When we write $v_{g_k}=\frac{1}{r_{g_k}}(a_{1,k},\dots,a_{n,k})$ (see (\ref{age})), we define the degree of $x_j$ as
\begin{equation}\label{deg}
\deg(x_j)=(a_{j,1},\dots,a_{j,1})\in\mathbb{Z}^m.
\end{equation}
Let us consider the Laurent polynomial ring $S_G=\C[V][t_1^{\pm1},\dots,t_m^{\pm1}]$ over $\C[V]$ with the $\mathbb{Z}^m$-grading given by setting $\deg(t_i)=1$ for all $i$. Let $\pi:X\to \C^n/G$ be an arbitrary relative minimal model. Then we can realize the Cox ring of $X$ as a graded subring of $S_G$. The following proposition follows from \cite[Lemma 3.3 and Proposition 3.4]{Y} for abelian groups $G\subset SL_n(\C)$.

\begin{prop}\label{cox}
With the notation above, $\mrm{Cox}(X)$ is isomorphic as a graded ring to the subring of $S_G$ generated by the homogeneous elements of the following subset
$$\{x_j t_1^{a_{j,1}}\cdots t_m^{a_{j,m}}\}_{j=1,\dots,n}\cup \{t_1^{-r_1},\dots,t_m^{-r_m}\}$$
of $S_G$ where $r_k=r_{g_k}$ is the order of $g_k$.
\end{prop}

\vspace{3mm}

\begin{rem}\label{Cl}
In the proposition the divisor class group $\mrm{Cl}(X(\Sigma))$ is identified with the subgroup of $\mathbb{Z}^m$ generated by $\{\deg(x_j)\}_{1\le j\le n}$ and $\{-r_k e_k\}_{1\le k\le m}$ where $e_1,\dots,e_m$ are the standard basis of $\mathbb{Z}^m$. Via the McKay correspondence for divisors (see Remark \ref{IR}), each junior element $g_k$ corresponds to a prime exceptional divisor $E_k$ of $\pi$, and the identification of $\mrm{Cl}(X(\Sigma))$ with the above subgroup of $\mathbb{Z}^m$ is chosen so that the divisor class of $E_k$ corresponds to $-r_k e_k\in\mathbb{Z}^m$. Note that the classes of the prime exceptional divisors $\{E_k\}$ give a $\Q$-basis of $\mathbb{Z}^m$ but not a basis unless $G$ is trivial.
\qed 
\end{rem}

\vspace{3mm}

Proposition \ref{cox} shows that $\mrm{Cox}(X)$ is a polynomial ring and thus $X$ is a toric variety satisfying $\mrm{Cox}(X)\cong S_X$. Note that $\mrm{Cox}(X)$ is independent of the choice of a relative minimal model $X$. Then Proposition \ref{quot} applies to $X=X(\Sigma)$ to show that there is an open subset $U_\Sigma\subset\mfk{X}:=\mrm{Spec}(\mrm{Cox}(X))$ which depends on $X$ (or $\Sigma$) such that $X$ is a geometric quotient of $U_\Sigma$ by the action of the torus $T:=\chi(\mrm{Cl}(X))$.

\section{Families of $G$-constellations over relative minimal models of $\C^n/G$}\label{4}

\subsection{Families of $G$-constellations}\label{4.1}

The aim of this paper is to try to realize a crepant resolution or more generally a relative minimal model $X\to \C^n/G$ as a (fine) moduli space of $G$-constellations. If this is actually possible, there must be a family $\mcl{F}$ of $G$-constellations on $X$ satisfying a certain universal property. As a candidate for such a family $\mcl{F}$, we consider a {\it geometrically natural family} of $G$-constellations (a {\it gnat family} for short) introduced by Logvinenko.

\begin{dfn}\label{def:gnat}(cf. \cite{L1},\cite{L3})
Let $\pi: X\to \C^n/G$ be a resolution and $q:\C^n\to \C^n/G$ the quotient map. A sheaf $\mcl{F}$ of $(R\rtimes G)\otimes_\C\mcl{O}_X$-module is called a {\bf gnat family} if the following conditions are satisfied:
\begin{itemize}
\item $\mcl{F}$ is locally free as an $\mcl{O}_X$-module, and its fiber $\mcl{F}_x$ is a $G$-constellation for any closed point $x\in X$
\item $q(\mrm{Supp}_{\C^n}\,\mcl{F}_x)=\pi(x)$ for any closed point $x\in X$ where $\mrm{Supp}_{\C^n}\,\mcl{F}_x\subset \C^n$ is the support of $\mcl{F}_x$ as an $\mcl{O}_{\C^n}$-module.\qed
\end{itemize}

\end{dfn}

Logvinenko showed that every gnat family is presented as a direct sum of invertible sheaves $\mcl{O}_X(-D_\chi)$ with $\Q$-divisors $\{D_\chi\}_{\chi\in \chi(G)}$ on $X$ parametrized by the characters of $G$ and that these $D_\chi$ must satisfy a certain condition called the reductor condition \cite[\S2]{L3}. In this paper we take exceptional $\Q$-divisors of a crepant resolution $\pi: X\to \C^n/G$ as $\{D_\chi\}_{\chi\in \chi(G)}$. Then each $D_\chi=D_{\rho_i}$ is written as $\sum_{k=1}^m b_{i,k}E_k$ with $E_k$ being the exceptional divisors of $\pi$ and $b_{i,k}\in\Q$. The reductor condition is translated into the condition that $b_{0,k}\in\Z$ for all $k=1,\dots,m$ and
\begin{equation}\label{gnat}
b_{i,k}+\frac{a_{j,k}}{r_k}-b_{k(i,j),k}\in\Z_{\ge0}
\end{equation}
for all $i\in\{0,\dots,r-1\},j\in\{1,\dots,n\}$ and $k$ where $a_{j,k}$ is the integer in (\ref{deg}). Since every gnat family is locally isomorphic to a normalized one i.e. $D_{\rho_0}=0$ as a $(R\rtimes G)\otimes_\C\mcl{O}_X$-module \cite[Corollary 3.5]{L3}, we always assume a gnat family is normalized from now on.

\begin{rem}
The set of gnat families is always nonempty \cite[\S3]{L3}. As we will see later, choosing a general point of $\Theta$ (avoiding finitely many hyperplanes) determines coefficients $\{b_{i,k}\}_i$ satisfying (\ref{gnat}) for each $k$ (when $X$ is a crepant resolution).
\qed
\end{rem}

Suppose that we are given a gnat family $\mcl{F}$ on a resolution $\pi: X\to \C^n/G$ or equivalently rational numbers $b_{i,k}$ satisfying the condition (\ref{gnat}). In this subsection we describe $G$-constellations in $\mcl{F}$ as quiver representations. Recall that a $G$-constellation is determined by assigning complex numbers $p_{i,j}$ to arrows $\alpha_{i,j}$ of the McKay quiver $Q$ so that the relations $p_{k(i,j),j'}p_{i,j}-p_{k(i,j'),j}p_{i,j'}=0$ are satisfied (see Section \ref{2}).

From now on we assume $X\to \C^n/G$ is crepant. Recall that $X$ admits $n+m$ homogeneous coordinates $\{x_j t_1^{a_{j,1}}\cdots t_m^{a_{j,m}}\}_{j=1,\dots,n}$ and $\{t_1^{-r_1},\dots,t_m^{-r_m}\}$ (see Subsection \ref{3.3}). If a point $x\in X$ is presented with these coordinates taking values 0 or 1, we call  such $x$ a {\it distinguished point}. In this case there is a unique cone $\sigma\in\Sigma$ such that the coordinates taking values 0 are associated with the rays of $\sigma$, and we say that a distinguished point $x$ is {\it associated} to $\sigma$. Such a point is denoted by $x_\sigma$.

The quiver representation $\mcl{F}_{x_\sigma}$ corresponding to the distinguished point $x_\sigma$ is obtained by setting $p_{i,j}=\delta_\sigma \epsilon_{\mcl{F},\sigma}$ with
\begin{equation}\label{delta}
\delta_\sigma=
\begin{cases}
0&\text{ if }\R_{\ge0} e_j\subset\sigma\\
1&\text{ otherwise}
\end{cases},\hspace{3mm}
\epsilon_{\mcl{F},\sigma}=
\begin{cases}
0&\text{ if }b_{i,k}+\frac{a_{j,k}}{r_k}-b_{k(i,j),k}>0\text{ for some }k\text{ with }\tau_k\subset\sigma\\
1&\text{ otherwise}
\end{cases}
\end{equation}
where $\tau_k$ is the ray of $\Sigma$ corresponding to the exceptional divisor $E_k$.

The big torus $T^n\subset\C^n$ acts on the representation $\{p_{i,j}\}$ by $p_{i,j}\mapsto t_j p_{i,j}$ for $\mathbf{t}=(t_1,\dots,t_n)\in T^n$. Since every point of $X$ comes from a distinguished point $x_\sigma$ by the $T^n$-action, we can explicitly describe representations $\mcl{F}_x$ for all $x\in X$.

As a special case of distinguished points, we consider a $T^n$-fixed point of $X$. Such a point is obtained as the distinguished point of a maximal cone $\sigma\in\Sigma_\mrm{max}$. For each representation $\rho_i\in\mrm{Irr}(G)$, we can associate to it a unique Laurent monomial $m_i\in\C[M]$ whose divisor of zeros in $X_\sigma(\cong\C^n)$ is the same as the $\Q$-divisor $D_{\rho_i}\cap X_\sigma$. More specifically, we can choose $m_i$ so that $v_k(m_i)=b_{i,k}$ for each $k\in \sigma$ where $v_k:M\to \Q$ is the valuation defined by the divisor $E_k$. See \cite{L1} for explicit examples. Then the set $\{\chi^{m_i}\}_{\rho_i\in\mrm{Irr}(G)}$ of Laurent monomials becomes a $G${\it -brick} introduced in \cite{J} if $\mcl{F}_{x_\sigma}$ is indecomposable. However, $\mcl{F}_{x_\sigma}$ may be decomposable in general.

\vspace{3mm}

\begin{ex}\label{A3} (A gnat family admitting a decomposable $G$-constellation)

Let $G\subset SL_2(\C)$ be the cyclic group of order 4 generated by $g=\begin{pmatrix}i&0\\0&-i\end{pmatrix}$ where $i$ is the imaginary unit. Then $\C^2/G$ has $A_3$-singularity and is realized as a toric variety defined by the cone $\sigma_0=\mrm{Cone}(e_1,e_2)\subset (N_G)_\R$ (see Subsection \ref{3.2}). In this case we have $N_G=\Z^2+\Z\cdot  (\frac{3}{4},\frac{1}{4})$. We write $x=\chi^{m_1}$ and $y=\chi^{m_2}$ for the dual basis $m_1,m_2\in M$ of $e_1,e_2$.

The unique crepant (or minimal) resolution $X$ of $\C^2/G$ is obtained by subdividing $\sigma_0$ by adding the rays
$$\tau_1=\R_{\ge0} \left(\frac{3}{4},\frac{1}{4}\right),\,\tau_2=\R_{\ge0} \left(\frac{1}{2},\frac{1}{2}\right),\,\tau_3=\R_{\ge0} \left(\frac{1}{4},\frac{3}{4}\right).$$
Let $E_k\subset X$ be the exceptional divisor corresponding to $\tau_k$.

We define a gnat family $\mcl{F}$ on $X$ by setting
$$\begin{aligned}
&D_{\rho_0}=0,&  &D_{\rho_1}=\frac{3}{4}E_1+\frac{1}{2}E_2-\frac{3}{4}E_3,\\
&D_{\rho_2}=\frac{1}{2}E_1+E_2-\frac{1}{2}E_3,& &D_{\rho_3}=\frac{1}{4}E_1+\frac{1}{2}E_2-\frac{1}{4}E_3
\end{aligned}$$
where $\rho_i$ is the $i$-th power of the character $\chi_1:G\to\C^* ;g\mapsto i$ defined by $x$. Then the Laurent monomials for the cone $\sigma_1:=\mrm{Cone}(\tau_1,\tau_2)$ is $\{1,x,y^2,y\}$ and hence $\mcl{F}_{x_{\sigma_1}}$ is indecomposable. On the other hand, those for the cone $\sigma_2:=\mrm{Cone}(\tau_2,\tau_3)$ is $\{1,\frac{x^3}{y^2},\frac{x^4}{y^2},\frac{x^2}{y}\}$ and $\mcl{F}_{x_{\sigma_2}}$ decomposes into three irreducible summands.
\qed
\end{ex}

\vspace{3mm}

\begin{rem}\label{nonproj2}
Theorem \ref{main} asserts that the indecomposability of the $G$-constellations $\mcl{F}_x$ is what we need for $X$ to be a fine moduli space. The example of \cite[\S5]{L2} gives a non-projective example of $X$ with $\mcl{F}$ parametrizing indecomposable $G$-constellations. Note that we only have to consider the torus-fixed points of $X$ in order to check the indecomposability of $\mcl{F}_x$ for all $x\in X$.
\qed
\end{rem}

\vspace{3mm}

So far we have considered $G$-constellation families over a resolution $X\to \C^n/G$. We can extend the notion of a gnat family to one over partial resolutions. In view of the moduli theory (particularly non-fine moduli spaces), it seems more natural to weaken Definition \ref{def:gnat} so that we allow non locally free sheaves when we consider singular varieties. Indeed, for a relative minimal model $X\to \C^n/G$, we can define a gnat family $\mcl{F}:=\bigoplus_{\rho_i\in\mrm{Irr}(G)}\mcl{O}_X(-D_{\rho_i})$ in this broader sense by choosing coefficients $b_{i,k}$ satisfying (\ref{gnat}) for each exceptional divisor $E_k$, similarly to the smooth case. Note that $\mcl{O}_X(-D_{\rho_i})$ in this case is a divisorial sheaf which is invertible on the smooth part of $X$ but possibly not on the whole $X$.  We also call such $\mcl{F}$ a {\it gnat family} on $X$ from now on. Since any relative minimal model $X$ is toric as well as a crepant resolution, we can explicitly calculate the $G$-constellations for singular points from $b_{i,k}$'s by using (\ref{delta}).

\vspace{3mm}

\begin{ex}\label{singular}(A gnat family on a singular relative minimal model)

A simplest example of a singular relative minimal model is given by the group $G=\{\pm I_4\}\subset SL_4(\C)$. In this case $\C^4/G$ itself is the relative minimal model since $G$ has no junior elements. Since it has no exceptional divisors, it has a uniquely determined gnat family $\mcl{F}$, which is isomorphic to $\mcl{O}_X\oplus \mcl{L}$ as a coherent sheaf where $\mcl{L}$ is the divisorial sheaf associated to the nontrivial element of $\mrm{Cl}(X)\cong\Z/2\Z$. Then $\mcl{F}$ parametrizes indecomposable $G$-constellations on the smooth part of $X$ while the fiber $\mcl{F}_O$ at the isolated singular point $O\in X$ decomposes into two 1-dimensional representations.

The space $\Theta$ of stability conditions is 1-dimensional, and one can check that the moduli spaces $\mcl{M}_\theta$ (including the $G$-Hilbert scheme) for generic $\theta$ are isomorphic to the blow-up of $X=\C^4/G$ at the singular point $O$. Therefore, we cannot realize $X$ as $\bar{\mfk{M}}_\theta$ with generic $\theta$ in this case.
\qed
\end{ex}

\subsection{Cones in $\Theta$ associated to exceptional divisors}\label{4.2}

In this subsection we introduce and study certain cones in $\Theta$ associated to the divisors $E_1\dots,E_m$ and a gnat family $\mcl{F}$.

We first show that every gnat family parametrizes indecomposable $G$-constellations outside a codimension two locus using the crepantness of a relative minimal model. As explained in the previous subsection, a gnat family $\mcl{F}$ on a relative minimal model $X=X(\Sigma)$ of $\C^n/G$ is determined by choosing coefficients $b_{i,k}$ satisfying (\ref{gnat}) for each exceptional divisor. Let us fix an exceptional divisor $E_\ell$ and consider its distinguished point $x_{\tau_\ell}\in E_\ell$ where $\tau_\ell\in \Sigma$ is the ray corresponding to $E_\ell$. We simply write $\mcl{F}_\ell$ for the $G$-constellation at $x_{\tau_\ell}$.

\begin{lem}\label{crepant}
Suppose that we have a nontrivial decomposition $R=M\oplus M'$ of the regular representation of $G$ as a $G$-module such that $M$ is a subrepresentation of the $G$-constellation $\mcl{F}_\ell$ (i.e. $x_j\cdot M\subset M$ in $\mcl{F}_\ell$ for any $j$). Then for any $\rho_i\subset M'$ and $x_j$ with $\rho_{k(i,j)}\subset M$, the subspace $x_j\cdot \rho_i\subset \rho_{k(i,j)}$ is nonzero in $\mcl{F}_\ell$. In particular, $\mcl{F}_\ell$ is indecomposable.
\end{lem}

{\em Proof.} We assume $x_j\cdot \rho_i$ is zero in $\mcl{F}_\ell$ in order to deduce a contradiction.  Then we have $b_{i,\ell}+\frac{a_{j,\ell}}{r_k}-b_{k(i,j),\ell}>0$. By the crepantness of $E_\ell$, we have $\frac{1}{r_k}\sum_{j'=1}^n a_{j',k}=1$ and hence $b_{k(i,j),\ell}+\frac{1}{r_k}\sum_{j'\ne j}a_{j',\ell}=b_{i,\ell}$. This means that,  in the quiver for $\mcl{F}_\ell$ (namely, the quiver obtained by removing arrows $\alpha_{i,j}$ from $Q$ for each $(i,j)$ with $x_{i,j}=0$ in $\mcl{F}_\ell$), we can reach the vertex $v_{\rho_i}$ from the vertex $v_{\rho_{k(i,j)}}$ through the arrows corresponding to $x_{j'}$'s with $j'\ne j$. This is contrary to the assumption that $M$ is a subrepresentation.

The latter claim follows from the former. Indeed, the fact that the characters of $G$ defined by $x_1,\dots,x_n$ generate $\chi(G)$ implies the existence of $x_j$ with $x_j\cdot M'\not\subset M'$ whenever we are given a nontrivial decomposition $\mcl{F}_\ell=M\oplus M'$ as representations, which is a contradiction.
\qed

\vspace{3mm}

Assume that we have a nontrivial decomposition $R=M\oplus M'$ as a $G$-module with $M$ being a subrepresentation in $\mcl{F}_\ell$. Then we define a new gnat family
$\mfk{t}_{\ell,M}(\mcl{F})$ associated with the coefficients $b'_{i,k}$ defined by
$$b'_{i,k}=
\begin{cases}
b_{i,k}+1\hspace{5mm}&\text{ if }k=\ell,\,\rho_0\in M\text{ and }\rho_i\in M'\\
b_{i,k}-1\hspace{5mm}&\text{ if }k=\ell,\,\rho_0\in M'\text{ and }\rho_i\in M\\
b_{i,k}\hspace{5mm}&\text{ otherwise }
\end{cases}.
$$
One can check that $\mfk{t}_{\ell,M}(\mcl{F})$ again satisfies the condition (\ref{gnat}) using the crepantness of $X$ and that $M'$ is a subrepresentaion of $\mfk{t}_{\ell,M}(\mcl{F})$. In terms of quivers for general points of $E_\ell$, the operation $\mfk{t}_{\ell,M}$ corresponds to adding arrows $\alpha_{i,j}$ for all $i,j$ with $\rho_i\subset M$ and $\rho_{k(i,j)}\subset M'$, removing arrows $\alpha_{i,j}$ for all $i,j$ with $\rho_i\subset M'$ and $\rho_{k(i,j)}\subset M$, and making the remaining arrows unchanged. Therefore, this operation just switches sub and quotient representations.

In fact every gnat family is obtained by performing $\mfk{t}_{\ell,M}$ for various $\ell$ and $M$ repeatedly from any gnat family. We show this in a geometric way in terms of cones in $\Theta$. Let $\sigma$ be a cone of the fan $\Sigma$ associated to the toric variety $X$. We define $C_\mcl{F}(\sigma)\subset\Theta\cong\chi(\bar{T}^{r-1})_\R$ as the cone generated by the characters of $x_{i,j}$ such that $p_{i,j}\ne0$ for the $G$-constellation $\mcl{F}_{x_\sigma}=\{p_{i,j}\}_{i,j}$ (see Section \ref{2}). Note that $C_\mcl{F}(\sigma)$ is inside the half space $H_S^+:=\{\sum_{i\in S}\theta_i\ge0\}\subset\Theta$ for a subset $S\subsetneq\{0,1,\dots,r-1\}$ if and only if $\bigoplus_{i\in S}\rho_i\subset R$ is a subrepresentation of $\mcl{F}_{x_\sigma}$. In particular $C_\mcl{F}(\sigma)$ has the full dimension $r-1$ if and only if $\mcl{F}_{x_\sigma}$ is an indecomposable representation. More generally, the number of irreducible direct summands of $\mcl{F}_{x_\sigma}$ is equal to $r-\dim(C_\mcl{F}(\sigma))$.

\begin{lem}\label{lem:contain}
Let $\theta\in\Theta$ be any generic stability condition. Then starting from any gnat family, we obtain a gnat family $\mcl{F}$ such that $C_{\mcl{F}_k}$ contains $\theta$ for all $k$ by applying $\mfk{t}_{\ell,M}$ repeatedly for various $\ell$ and $M$.
\end{lem}

{\em Proof.} First note that $\Theta$ is divided into finitely many open cones by hyperplanes $H_S:=\{\sum_{i\in S}\theta_i=0\}\subset\Theta$ for subsets $S\subsetneq\{0,1,\dots,r-1\}$.

Take any gnat family $\mcl{F}$ and fix $\ell$. By Lemma \ref{crepant}, $C_{\mcl{F}_\ell}$ is of full dimension. Let $F$ be any facet of $C_{\mcl{F}_\ell}$. Then it must be defined as the intersection of $C_{\mcl{F}_\ell}$ and the hyperplane $H_S$ for some $S$. We may assume $C_{\mcl{F}_\ell}\subset H_S^+=\{\sum_{i\in S}\theta_i\ge0\}$ by replacing $S$ with its complement if necessary. We put $M:=\bigoplus_{i\in S}\rho_i$. Then $C_{\mfk{t}_{\ell,M}(\mcl{F})}$ also contains $F$ as a facet and is contained in $H_S^-=\{\sum_{i\in S}\theta_i\le0\}$. This follows from the fact that the operation $\mfk{t}_{\ell,M}$ does not change the arrows whose head and tail are both in $S$ or its complement in the quiver for $\mcl{F}_\ell$.

Since the number of chambers in $\Theta$ is finite, we can reach the chamber containing $\theta$ by repeatedly applying $\mfk{t}_{\ell,M}$ for various $M$. Then we just have to proceed the same things for all $\ell$.
\qed

\subsection{Construction of the map $\alpha:X\to \bar{\mfk{M}}_{X,\mcl{F}}$}\label{4.3}

In this subsection we will construct a candidate $\bar{\mfk{M}}_{X,\mcl{F}}$ of a moduli space parametrizing the fibers of a gnat family $\mcl{F}$ on a relative minimal model $X$ of $\C^n/G$. We will also construct a birational morphism $\alpha:X\to \bar{\mfk{M}}_{X,\mcl{F}}$ so that it is an isomorphism if $\bar{\mfk{M}}_{X,\mcl{F}}$ is actually a fine moduli space (cf. Theorem \ref{main}). For this, we first construct a homomorphism $\vphi_\mcl{F}:A/I\to \mrm{Cox}(X)$ of rings.

Recall that $A/I$ is generated by the variables $\bar{x}_{i,j}$, which correspond to arrows in the quiver $Q$ whose tails and heads are vertices for $\rho_i$ and $\rho_{k(i,j)}$ respectively (see Section \ref{2}). We define a ring homomorphism $\vphi'_\mcl{F}:A=\C[\{x_{i,j}\}]\to S_G=\C[V][t_1^{\pm1},\dots,t_m^{\pm1}]$ by
\begin{equation}\label{map}
\vphi'_\mcl{F}(x_{i,j})=x_j \prod_{k=1}^m t_k^{r_k(b_{k(i,j),k}-b_{i,k})}\in S_G
\end{equation}
where $r_k$ is the order of the junior element $g_k$ (see Subsection \ref{3.2}) and $b_{i,k}$ are the defining coefficients of $\mcl{F}$ (see Subsection \ref{4.1}). It is easy to see that $\vphi'_\mcl{F}$ is zero for elements in the defining ideal $I$. We can also show that $\vphi'_\mcl{F}$ factors through $\mrm{Cox}(X)$.

\begin{lem}
For any $i$ and $j$, $\vphi'_\mcl{F}(x_{i,j})$ belongs to $\mrm{Cox}(X)$. Therefore, $\vphi'_\mcl{F}$ induces a ring homomorphism $\vphi_\mcl{F}:A/I\to \mrm{Cox}(X)$.
\end{lem}

{\em Proof.} By the definitions of $k(i,j)$ and $b_{i,k}$, we see that $\chi_j(g_k)=e^{2\pi i (b_{k(i,j),k}-b_{i,k})}$ where $\chi_j$ is the character of $G$ defined by $x_j$ (cf. (\ref{seq2})). Thus, $r_k(b_{k(i,j),k}-b_{i,k})$ is equal to $a_{j,k}$ modulo $r_k\mathbb{Z}$ (see Section \ref{3}). We also have $a_{j,k}\ge0$ and $r_k(b_{k(i,j),k}-b_{i,k})\le a_{j,k}$. By considering the form of the generators of $\mrm{Cox}(X)$ in Proposition \ref{cox}, we see that $\vphi'(x_{i,j})$ is written as a product of these generators.
\qed

\vspace{3mm}

Recall that the divisor class group $\mrm{Cl}(X)$ is regarded as a subgroup of $\Z^n$ of finite index (see Remark \ref{Cl}). Thus, the torus $\bar{T}^m:=\chi(\mrm{Cl}(X))$ is a finite quotient of $T^m=(\C^*)^{\times m}$, which acts on $S_G$ so that each monomial $f\prod_{k=1}^m t_k^{a_k}\,(f\in \C[V])$ has weight $(a_1,\dots,a_k)$. We define a group homomorphism from $\bar{T}^m$ to the acting torus $\bar{T}^{r-1}=(\C^*)^{\times r}/\Delta$ on $A/I$ (see Section \ref{2}) by
$$\psi_\mcl{F}:\bar{T}^m\to \bar{T}^{r-1};\,(\overline{t_1,\dots,t_m})\mapsto \left(\overline{\prod_{k=1}^m t_k^{r_k b_{0,k}},\dots,\prod_{k=1}^m t_k^{r_k b_{r-1,k}}}\right).$$
Note that the two tori naturally act also on the spectra $\mfk{X}=\mrm{Spec}\,\mrm{Cox}(X)$ and $\mcl{R}=\mrm{Spec}\,A/I$ respectively.

\begin{lem}\label{equi}
The induced morphism $\bar{\vphi}_\mcl{F}:\mfk{X}\to\mcl{R}$ from $\vphi_\mcl{F}:A/I\to \mrm{Cox}(X)$ is equivariant with respect to the group homomorphism $\psi_\mcl{F}:\bar{T}^m\to \bar{T}^{r-1}$.
\end{lem}

{\em Proof.} The desired equivariance of $\bar{\vphi}_\mcl{F}$ is equivalent to
$$\vphi_\mcl{F}(\psi_\mcl{F}(\mathbf{t})\cdot f)=\mathbf{t}\cdot \vphi_\mcl{F}(f)$$
for any $\mathbf{t}\in \bar{T}^m$ and $f\in A/I$. It suffices to check this for the generators $\bar{x}_{i,j}\in A/I$. This can be done by (\ref{map}) noticing that $(t_0,\dots,t_{r-1})\in T^r$ acts on $\bar{x}_{i,j}$ as multiplication of $t_i^{-1}t_{k(i,j)}$ (see (\ref{Tr})).
\qed

\vspace{3mm}

\begin{rem}
The construction of $\bar{\vphi}_\mcl{F}$ and $\psi_\mcl{F}$ depends only on (the coefficients of) $\mcl{F}$ but not on the choice of a relative minimal model $X$.
\qed
\end{rem}

\vspace{3mm}

\begin{ex}
The homomorphism $\vphi_\mcl{F}:\mcl{R}\to \mrm{Cox}(X)=\C[xt_1^3 t_2^2 t_3,yt_1 t_2^2 t_3^3,t_1^{-4},t_2^{-4},t_3^{-4}]$ for the case in Example \ref{A3} is explicitly computed as follows:
$$\begin{aligned}
&\vphi_\mcl{F}(x_0)=xt_1^3 t_2^2 t_3^{-3},& &\vphi_\mcl{F}(x_1)=xt_1^{-1} t_2^2 t_3,& &\vphi_\mcl{F}(x_2)=xt_1^{-1} t_2^{-2} t_3,& &\vphi_\mcl{F}(x_3)=xt_1^{-1} t_2^{-2} t_3,\\
&\vphi_\mcl{F}(y_0)=yt_1 t_2^2 t_3^{-1},& &\vphi_\mcl{F}(y_1)=yt_1^{-3} t_2^{-2} t_3^3,& &\vphi_\mcl{F}(y_2)=yt_1 t_2^{-2} t_3^1,& &\vphi_\mcl{F}(y_3)=yt_1 t_2^2 t_3^{-1}
\end{aligned}$$
where we write $x_j=\bar{x}_{1,j}$ and $y_j=\bar{x}_{2,j}$ for the readability. Moreover, if we identify $\chi(\bar{T}^{r-1})$ with $\Z^3$ by
$$\deg(x_0)\mapsto (1,0,0),\deg(x_1)\mapsto (-1,1,0),\deg(x_2)\mapsto (0,-1,1),\deg(x_3)\mapsto (0,0,-1)$$
and also identify $\chi(\bar{T}^m)$ with the subgroup of $\Z^3$ generated by $\deg(f)$ for homogeneous $f\in\mrm{Cox}(X)$ (see Remark \ref{Cl}), then the group homomorphism $\psi_\mcl{F}^*:\chi(\bar{T}^{r-1})\to \chi(\bar{T}^m)$ becomes a linear map
$$\Z^3\to \Z^3;\;
\mathbf{v}\mapsto
\begin{pmatrix}
3&2&1\\
2&4&2\\
-3&-2&-1
\end{pmatrix}
\mathbf{v}.$$
So $(\psi_\mcl{F}^*)_\R$ may not be surjective in general.
\qed
\end{ex}

\vspace{3mm}

Since $\mfk{X}$ is irreducible, the morphism $\bar{\vphi}_\mcl{F}$ factors through an irreducible component of $\mcl{R}$. Recall that there is a unique irreducible component $\mcl{V}$ of $\mcl{R}(\subset \C^{nr})$ which does not lie in any coordinate hyperplane of $\C^{nr}$. Since $\bar{\vphi}_\mcl{F}(1,\dots ,1)\in \mcl{V}$, we have $\mrm{Im}\,\bar{\vphi}_\mcl{F}\subset \mcl{V}$.

We will introduce open subsets of $\mfk{X}$ and $\mcl{R}$ so that the restriction of $\bar{\vphi}_\mcl{F}$ to these open subsets induce a morphism between $X$ and a candidate of a moduli space parametrizing the fibers of a gnat family $\mcl{F}$. We fix a relative minimal model $X$ of $\C^n/G$ from now on. As mentioned in the previous section, $X$ is a toric variety and is determined by a fan $\Sigma\subset (N_G)_\mathbb{R}$. Then, by Proposition \ref{fan}, each $\sigma\in\Sigma_\mrm{max}$ is generated by $n$ elements in the set $C=\{e_1,\dots,e_n,v_{g_1},\dots,v_{g_m}\}\subset N_G$ where $\{e_j\}$ is the standard basis of $N\subset N_G$ and $v_{g_k}$ is the vector associated to the junior element $g_k\in G$. Let $y_j=x_j \prod_{k=1}^m t_k^{a_{j,k}}\,(j=1,\dots,n)$ and $z_k=t_k^{-r_k}\,(k=1,\dots,m)$ be the generators of $\mrm{Cox}(X)$ in Proposition \ref{cox}. We say that $y_j$ (resp. $z_k$) is the coordinate of $\mfk{X}$ which is {\em assigned } to $e_j$ (resp. $v_{g_k}$). For each $\sigma\in\Sigma_\mrm{max}$, take the unique generators of it from $C$, and let $\mfk{X}_\sigma\subset \mfk{X}$ be the open subset defined as the non-vanishing locus of the coordinates assigned to these generators. Then the open subset $\mfk{X}_\Sigma\subset\mfk{X}$ defined as
$$\mfk{X}_\Sigma=\bigcup_{\sigma\in \Sigma_\mrm{max}}\mfk{X}_\sigma\subset \mfk{X}$$
is the same as the open subset $U_{\tilde{\Sigma}}$ in Subsection \ref{3.3} associated to the toric variety $X(\Sigma)$.Thus, $X=X(\Sigma)$ is obtained as the geometric quotient of $\mfk{X}_\Sigma$ by the action of $\bar{T}^m$.

We also define an open subset $\mcl{V}_{\Sigma,\mcl{F}}\subset \mcl{V}$ as follows. For each coordinate $\bar{x}_{i,j}$ of $\mcl{R}$, we can write $\vphi(\bar{x}_{i,j})$ uniquely as a product of coordinates $y_j$ and $z_k$ of $\mfk{X}$. For each $\sigma\in\Sigma_\mrm{max}$, we define an open subset $\mcl{V}_{\sigma,\mcl{F}}\subset \mcl{V}$ as the non-vanishing locus of $\bar{x}_{i,j}$'s for which all the factors of $\vphi_\mcl{F}(\bar{x}_{i,j})$ are coordinates assigned to generators of $\sigma$:
$$\mcl{V}_{\sigma,\mcl{F}}=\{v\in\mcl{V}\mid \bar{x}_{i,j}(v)\ne0\text{ for any }i,j\text{ such that }\vphi_\mcl{F}(\bar{x}_{i,j})\text{ is invertible on }X_\sigma\}.$$

Then we define an open subset $\mcl{V}_{\Sigma,\mcl{F}}\subset\mcl{V}$ as
$$\mcl{V}_{\Sigma,\mcl{F}}=\bigcup_{\sigma\in \Sigma_\mrm{max}}\mcl{V}_{\sigma,\mcl{F}}.$$

In \cite[Theorem 3.10]{CMT} it is shown that $\mcl{V}$ is a possibly nonnormal toric variety, namely $\mcl{V}$ contains an algebraic torus $T_\mcl{V}$ as an open dense subset and the natural action of $T_\mcl{V}$ on itself extends to an action on the whole $\mcl{V}$. If $\mcl{V}$ is normal, it is a usual toric variety treated in Subsection \ref{3.1}. In fact the action of $T_\mcl{V}$ on $\mcl{V}$ is realized as a combination of actions of the two tori $\bar{T}^{r-1}$ and $T^n$ where $T^n=(\C^*)^{\times n}\subset \C^n$ is the standard torus \cite[Remark 3.12]{CMT}. Note that $T^n$ naturally acts on both $\mcl{R}$ and $\mcl{V}$ and that $\bar{\vphi}_\mcl{F}$ is clearly equivariant with respect to these actions. This implies that $\mcl{V}_{\Sigma,\mcl{F}}$ is $T_\mcl{V}$-invariant and thus $\mcl{V}_{\Sigma,\mcl{F}}$ is also a possibly nonnormal toric variety.

Next we consider the fan of the normalization of $\mcl{V}_{\Sigma,\mcl{F}}$. Let $N_\mcl{V}=\mrm{Hom}_\mrm{alg.gp.}(\C^*,T_\mcl{V})$ be the free abelian group of one-parameter subgroups of $T_\mcl{V}$. The functions $\bar{x}_{i,j}$ on $\mcl{V}$ are regarded as generators of the dual abelian group $M_\mcl{V}:=\mrm{Hom}_\Z(N_\mcl{V},\Z)$. For any $\sigma\in\Sigma$, let $S_\sigma\subset M_\mcl{V}$ be the semigroup consisting of the homogeneous regular functions on $\mcl{V}_{\sigma,\mcl{F}}$, that is, $S_\sigma$ is generated by $\bar{x}_{i,j}$ for all $i,j$ and by $\bar{x}_{i,j}^{-1}$ for which $\vphi_\mcl{F}(\bar{x}_{i,j})$ is invertible on $X_\sigma$. We define a cone $\tilde{\sigma}$ in $(N_\mcl{V})_\R=N_\mcl{V}\otimes \R$ as the dual of $S_\sigma$:
$$\tilde{\sigma}=\{u\in (N_\mcl{V})_\R\mid f(u)\ge0,\forall f\in S_\sigma\}.$$

\begin{lem}\label{fan1}
The maximal dimensional cones (whose dimensions may be different) of the fan $\tilde{\Sigma}$ associated to the normalization $\bar{\mcl{V}}_{\Sigma,\mcl{F}}$ of $\mcl{V}_{\Sigma,\mcl{F}}$ are the cones in the set $\{\tilde{\sigma}\}_{\sigma\in\Sigma_{\mrm{max}}}$ which are maximal with respect to inclusion. 
\end{lem}

{\em Proof.} Note that the normalization $\bar{\mcl{V}}_{\sigma,\mcl{F}}$ of the affine variety $\mcl{V}_{\sigma,\mcl{F}}$ is obtained as the affine toric variety $\mrm{Spec}\,\C[\tilde{\sigma}^\vee\cap M_\mcl{V}]$ where $\tilde{\sigma}^\vee=\{f\in (M_\mcl{V})_\R\mid f(u)\ge0,\forall u\in \tilde{\sigma}\}$. Since $\bar{\mcl{V}}_{\Sigma,\mcl{F}}$ is covered by $\bar{\mcl{V}}_{\sigma,\mcl{F}}$'s, the maximal dimensional cones of $\tilde{\Sigma}$ are exhausted by the maximal cones in $\{\tilde{\sigma}\}_{\sigma\in\Sigma}$.
\qed

\vspace{3mm}

\begin{rem}
The set $\{\tilde{\sigma}\}_{\sigma\in\Sigma_{\mrm{max}}}$ contains non-maximal cones in general (see Remark \ref{example}). If all $\bar{T}^{r-1}$-orbits in $\mcl{V}_{\Sigma,\mcl{F}}$ (or $\bar{\mcl{V}}_{\Sigma,\mcl{F}}$) are closed, this does not happen (cf. Lemma \ref{closed}).
\qed
\end{rem}

\vspace{3mm}

Next we construct a variety $\bar{\mfk{M}}_{X,\mcl{F}}$ which is a candidate of a (fine) moduli space parameterizing of $G$-constellations appearing as fibers $\mcl{F}_x$. We do this by taking a quotient of $\bar{\mcl{V}}_{\Sigma,\mcl{F}}$ by the $\bar{T}^{r-1}$-action. However, we cannot construct $\bar{\mfk{M}}_{X,\mcl{F}}$ as the orbit space of $\bar{\mcl{V}}_{\Sigma,\mcl{F}}$ in general since $\bar{\mcl{V}}_{\Sigma,\mcl{F}}$ may have non-closed $\bar{T}^{r-1}$-orbits. We shall show that $\bar{\mfk{M}}_{X,\mcl{F}}$ still exists as a categorical quotient in the category of algebraic varieties. We also show that $\bar{\mfk{M}}_{X,\mcl{F}}$ is a toric variety and describe its fan.

In general a good quotient of a normal toric variety by the action of a subtorus of the big torus does not exist as a variety. However, a categorical quotient always exists in the category of normal toric varieties and toric morphisms \cite{AH1}. Note that the abelian group $N_G$ is isomorphic to the quotient of $N_\mcl{V}$ by the subspace consisting of one-parameter subgroups of $\bar{T}^{r-1}\subset T_\mcl{V}$, which follows from the fact that the $\bar{T}^{r-1}$-invariant part of $M_\mcl{V}$ is the same as $M_G$. Let $q_G:(N_\mcl{V})_\R\to (N_G)_\R$ be the projection. Then the fan $\Sigma_\mcl{F}\subset(N_G)_\R$ corresponding to the toric quotient is explicitly constructed as follows (cf. Proof of \cite[Theorem 2.3]{AH1}):

(1) Let $S$ be the set of cones $q_G(\tilde{\sigma}),\,\tilde{\sigma}\in \tilde{\Sigma}$ which are maximal with respect to inclusion.

(2) As long as there exist cones $\sigma_1$ and $\sigma_2$ in $S$ such that $\sigma_1\cap \sigma_2$ is not a face of $\sigma_1$, do the following: Let $\tau_2$ be the minimal face of $\sigma_2$ that contains $\sigma_1\cap \sigma_2$. If $\tau_2 \not\subset\sigma_1$, replace $\sigma_1$ by the cone $\mrm{Cone}(\sigma_1, \tau_2)$ generated by $\sigma_1$ and $\tau_2$. Otherwise let $\tau_1$ be the minimal face of $\sigma_1$ that contains $\sigma_1\cap \sigma_2$ and replace $\sigma_2$ by $\mrm{Cone}(\sigma_2, \tau_1)$. Omit all cones of $S$ which are properly contained in the new one.

(3) If the process (2) cannot be carried out anymore, we obtain $\Sigma_\mcl{F}$ as the fan consisting of all the faces of the cones of $S$.

The fan $\Sigma_\mcl{F}$ has the following property: For any maximal cone $\sigma\in\Sigma_\mcl{F}$, the intersection $q_G^{-1}(\sigma)\cap |\tilde{\Sigma}|$ is a union of cones in $\tilde{\Sigma}$. $\Sigma_\mcl{F}$ is characterized as the finest fan among the set of fans in $(N_G)_\R$ satisfying this property.

\begin{dfn}\label{candidate}
We define $\bar{\mfk{M}}_{X,\mcl{F}}$ as the toric variety corresponding to the fan $\Sigma_\mcl{F}\subset(N_G)_\R$ constructed as above.
\end{dfn}

To summarize, the toric morphism $q_2:\bar{\mcl{V}}_{\Sigma,\mcl{F}}\to\bar{\mfk{M}}_{X,\mcl{F}}$ induced by $q_G$ is the categorical quotient by the $\bar{T}^{r-1}$-action in the toric category.

\vspace{3mm}

\begin{rem}
The categorical quotient $q_2:\bar{\mcl{V}}_{\Sigma,\mcl{F}}\to\bar{\mfk{M}}_{X,\mcl{F}}$ in the toric category may not be a good quotient since it is possibly not an affine morphism. If $q_G^{-1}(\sigma)\cap |\tilde{\Sigma}|$ is a single cone of $\tilde{\Sigma}$ for every $\sigma\in\Sigma_\mcl{F}$, then $q_2$ is a good quotient (cf. \cite[Proposition 3.2]{AH1}).
\qed
\end{rem}

\vspace{3mm}

\begin{prop}\label{alpha}
The $\bar{T}^{r-1}$-invariant morphism $q_2$ is the categorical quotient in the category of algebraic varieties. Moreover, the morphism $\bar{\vphi}_\mcl{F}:\mfk{X}\to\mcl{R}$ induces a proper birational toric morphism $\alpha:X\to\bar{\mfk{M}}_{X,\mcl{F}}$ which makes the following diagram commutative:\\
\begin{equation}\label{morph}
\begin{CD}
\mfk{X}_{\Sigma} @>\bar{\vphi}_{\Sigma,\mcl{F}}>> \bar{\mcl{V}}_{\Sigma,\mcl{F}}\\ @V q_1 VV @VV q_2 V\\ X@>\alpha>> \bar{\mfk{M}}_{X,\mcl{F}}\\ 
\end{CD}
\end{equation}

\end{prop}

{\it Proof.} We first show that $q_2$ is the categorical quotient. For any $\sigma\in\Sigma_\mcl{F}$, the open subset $\tilde{U}_\sigma:=q_2^{-1}(\bar{U}_\sigma)$ is torus-invariant and thus is associated with a subfan $\Sigma'$ of $\tilde{\Sigma}$. Let $U'\subset \bar{\mcl{V}}_{\Sigma,\mcl{F}}$ be the open subset corresponding to the cone $\sigma'$ generated by the cones in $\Sigma'$. Note that we have $H^0(\mcl{O}_{U'})=H^0(\mcl{O}_{\tilde{U}_\sigma})$ and $q_G(\sigma')=\sigma$. Thus we also have $H^0(\mcl{O}_{\bar{U}_\sigma})=H^0(\mcl{O}_{U'})^{\bar{T}^{r-1}}=H^0(\mcl{O}_{\tilde{U}_\sigma})^{\bar{T}^{r-1}}$. By \cite[Corollary 1.4]{ACH}, $\bar{U}_\sigma$ is the categorical quotient of $\tilde{U}_\sigma$ for the $\bar{T}^{r-1}$-action in the category of algebraic varieties. Note that $q_2$ is surjective since, for each cone $\sigma''\in \Sigma$ with $\sigma''\subset \sigma$, the associated cone $\tilde{\sigma}''\in\tilde{\Sigma}$ satisfies $q_G(\tilde{\sigma}'')\supset \sigma''$. Thus the categorical quotient of $\tilde{U}_\sigma$ is strong in the sense that, for every open $U_0\subset \bar{U}_\sigma$, the restriction $q_2^{-1}(U_0)\to U_0$ is a categorical quotient. This implies that $q_2:U\to \bar{\mfk{M}}_{X,\mcl{F}}$ is globally a categorical quotient.

The last claim follows from the fact that the following commutative diagram of $\R$-vector spaces induces a commutative diagram of the fans corresponding to the four varieties in the diagram (\ref{morph}).
\begin{equation}\label{diag1PS}
\begin{CD}
(N_{\mfk{X}})_\R @>\bar{\psi}_{\mcl{F}}>> (N_{\mcl{V}})_\R \\
@V p_G VV @VV q_G V\\
(N_G)_\R @= (N_G)_\R.\\ 
\end{CD}
\end{equation}
where $N_{T_\mfk{X}}$ is the rank-$(m+n)$ lattice associated to the affine space $\mfk{X}$ and the map $\bar{\psi}_{\mcl{F}}$ is induced by $\psi_\mcl{F}:\bar{T}^m\to \bar{T}^{r-1}$.

Then the morphism $\alpha$ corresponds to a refinement $\Sigma\to\Sigma_\mcl{F}$, and thus it is proper and birational. 
\qed

\vspace{3mm}

\begin{rem}
1. Although it is more natural to construct $\bar{\mfk{M}}_{X,\mcl{F}}$ by taking a categorical quotient $\mfk{M}_{X,\mcl{F}}$ of $\mcl{V}_{\Sigma,\mcl{F}}$ and then taking its normalization, we avoid this since constructing a categorical quotient is harder in nonnormal cases. 
\\
2. As alternatives of  the categorical quotient in the category of algebraic varieties, one may consider quotients in categories of not-necessarily-separated (toric) schemes. There are studies of quotients in such enlarged categories, particularly for toric varieties (see e.g. \cite{AH2}). We do not go further in this direction since we are mainly interested in cases where $\bar{\mfk{M}}_{X,\mcl{F}}$ is a fine moduli space (and in particular $\bar{\mcl{V}}_{\Sigma,\mcl{F}}$ admits a geometric quotient), and in such cases these notions of quotients coincide.
\qed
\end{rem}

\vspace{3mm}

The following lemma shows that any gnat family on $X$ minus a codimension-two locus parametrizes mutually non-isomorphic $G$-constellations. 

\begin{lem}\label{crepant divisor}
Let $\mcl{F}$ be a gnat family on a relative minimal model $X=X(\Sigma)$ of $\C^n/G$ and let $\mcl{F}_\ell$ be the $G$-constellation for a general point of an irreducible exceptional divisor $E_\ell\subset X$. Then the closure $Y$ of the $T_\mcl{V}$-orbit of $\mcl{F}_\ell$ (as a subset of $\mcl{V}_{\Sigma,\mcl{F}}$) is a torus-invariant irreducible divisor. In particular it defines a ray $\tau_{\mcl{F},\ell}$ of the fan of the normal toric variety $\bar{\mcl{V}}_{\Sigma,\mcl{F}}$.
\end{lem}

{\em Proof.} Writing $\mcl{F}_\ell$ as $\{p_{i,j}\}_{i,j}$, $Y$ is regarded as the zero locus of $x_{i,j}$'s in $\mcl{V}$ for $i,j$ with $p_{i,j}=0$. If a coordinate $x_{i,j}$ is zero, then the commutative relation $x_{k(i,j),j'}x_{i,j}-x_{k(i,j'),j}x_{i,j'}=0$ implies $x_{k(i,j'),j}x_{i,j'}=0$. By the description of $\mcl{F}_\ell$ explained in the proof of  Lemma \ref{crepant}, we see that exactly one, say $x'$, of the two coordinates $x_{k(i,j'),j}$ and $x_{i,j'}$ is zero on $Y$. This means that the zero locus $\{x_{i,j}=x'=0\}\subset \mcl{V}$ is of pure codimension one. Since the quiver for $\mcl{F}_\ell$ is connected by Lemma \ref{crepant} again, we can choose commutative relations of the coordinates successively so that the zero locus of all $x_{i,j}$'s for which $p_{i,j}=0$ is also of pure codimension one. Since $Y$ is the closure of a single $T_\mcl{V}$-orbit, it is irreducible and the claim follows. 
\qed

\vspace{3mm}

Since the cones $C_{\mcl{F}_\ell}$ with various families $\mcl{F}$ cover $\Theta$ for each $\ell$ by Lemma \ref{lem:contain}, we obtain the following corollary: 

\begin{cor}\label{cor}
Let $\theta\in\Theta$ be a generic stability condition. Then the normalized fine moduli space $\bar{\mfk{M}}_\theta$ contains all crepant divisors $E_1,\dots,E_m$. In other words, the fan of $\bar{\mfk{M}}_\theta$ contains rays generated by $v_g\in N_G$ for all junior elements $g\in G$ (see Definition \ref{junior}).
\end{cor}

\section{Proof of the main result}\label{5}

Let $X\to \C^n/G$ be a crepant resolution for a finite abelian group $G\subset SL_n(\C)$ and let $\mcl{F}$ be a gnat family on $X$. In this subsection we study when the birational morphism $\alpha:X\to\bar{\mfk{M}}_{X,\mcl{F}}$ constructed in the previous section is an isomorphism to a fine moduli space.

Firstly we should note that the categorical quotient $\bar{\mfk{M}}_{X,\mcl{F}}$ may not be a geometric (nor even good) quotient, that is, $\bar{\mcl{V}}_{\Sigma,\mcl{F}}$ (or $\bar{\mcl{V}}_{\Sigma,\mcl{F}}$) contains non-closed $\bar{T}^{r-1}$-orbits in general. We can interpret the (non-)closedness of a $\bar{T}^{r-1}$-orbit in terms of corresponding $G$-constellations as follows.

\begin{lem}\label{closed}
The following four conditions are equivalent:\\
(a): Every $\bar{T}^{r-1}$-orbit in $\mcl{V}_{\Sigma,\mcl{F}}$ is closed.\\
(b): Every $\bar{T}^{r-1}$-orbit in $\bar{\mcl{V}}_{\Sigma,\mcl{F}}$ is closed.\\
(c): For each $p\in \mcl{V}_{\Sigma,\mcl{F}}$, the stabilizer subgroup $\mrm{Stab}_p(\bar{T}^{r-1})\subset \bar{T}^{r-1}$ is finite.\\
(d): For each $p\in \mcl{V}_{\Sigma,\mcl{F}}$, the corresponding $G$-constellation $M_p$ (see Section\ref{2}) is indecomposable.
\end{lem}

{\em Proof.} (b)$\Rightarrow$(c): By assumption (a), the quotient map $q:\bar{\mcl{V}}_{\Sigma,\mcl{F}}\to \bar{\mfk{M}}_{X,\mcl{F}}$ is geometric and thus $q^{-1}q(p)$ is isomorphic to the homogeneous space $\bar{T}^{r-1}/\mrm{Stab}_p(\bar{T}^{r-1})$ for any point $p\in \bar{\mcl{V}}_{\Sigma,\mcl{F}}$. Since $\mrm{Stab}_{p_0}(T_r)$ are trivial for points $p_0$ in the big torus of $\mcl{V}_{\Sigma,\mcl{F}}$, we have $\dim (\bar{T}^{r-1}/\mrm{Stab}_p(T_r))\ge\dim \bar{T}^{r-1}$ for all $p$ by semicontinuity and thus $\mrm{Stab}_p(T_r)$ must be finite.

(c)$\Rightarrow$(d): Assume $M_p$ admits a nontrivial decomposition $M_p\cong M_1\oplus M_2$ as a representation of the quiver $Q$ (or as a $G$-equivariant $\C[V]$-module). This particularly gives a nontrivial decomposition as a $G$-module, and accordingly we consider the decomposition $(\C^*)^r=T_1\times T_2$ of the torus (i.e. $T_i=\prod_{\rho\subset M_i} \C^*$). Then the 1-parameter subgroup $\lambda:\C^*\hookrightarrow T_1\times T_2;\,t\mapsto((1,\dots,1),(t,\dots,t))$ acts trivially on $M_p$ and satisfies $\lambda(\C^*)\cap \Delta=\{1\}$. Therefore, $\mrm{Stab}_p(\bar{T}^{r-1})$ has positive dimension.

(d)$\Rightarrow$(a): Assume that the orbit $\bar{T}^{r-1}\cdot p_1$ of a point $p_1$ is not closed in $\mcl{V}_{\Sigma,\mcl{F}}$, that is, there exists $p_2\in\mcl{V}_{\Sigma,\mcl{F}}$ such that $p_2\in \overline{\bar{T}^{r-1}\cdot x}\setminus \bar{T}^{r-1}\cdot x$. Then there is a nontrivial filtration
$$ F_\bullet: F_0=0\subsetneq F_1 \subsetneq \cdots \subsetneq F_{\ell-1}\subsetneq F_\ell=M_{p_1}\;(\ell\ge2)$$
of subrepresentations of the $G$-constellation $M_{p_1}$ such that the $G$-constellation $M_{p_2}$ is isomorphic to $\mrm{Gr}(F_\bullet):=\bigoplus_{i=1}^\ell F_i/F_{i-1}$. This follows from the same argument as \cite[Theorem 2.8]{Kir}. In particular, $M_{p_2}$ admits a nontrivial decomposition as a $G$-equivariant $\C[V]$-module.

(a)$\iff$(b): This follows from the fact that the normalization map $\bar{\mcl{V}}_{\Sigma,\mcl{F}}\to\mcl{V}_{\Sigma,\mcl{F}}$ naturally gives one-to-one correspondence between the sets of $\bar{T}^{r-1}$-orbits of $\bar{\mcl{V}}_{\Sigma,\mcl{F}}$ and $\mcl{V}_{\Sigma,\mcl{F}}$ \cite[Theorem 3.A.3]{Co2}.
\qed

\vspace{3mm}

\begin{rem}\label{example}
In Example \ref{A3} we saw a gnat family $\mcl{F}$ admitting a decomposable fiber $\mcl{F}_{x_{\sigma_2}}$. One can check that the fibers $\mcl{F}_{x_{\tau_k}},\,k=2,3$ admit filtrations such that their associated graded representation is equal to $\mcl{F}_{x_{\sigma_2}}$. Thus, the $\bar{T}^{r-1}$-orbits of $\mcl{F}_{x_{\tau_k}}$ are non-closed in $\mcl{V}_{\Sigma,\mcl{F}}$ and their closures contain the orbit of $\mcl{F}_{x_{\sigma_2}}$. This implies that the divisors $E_2$ and $E_3$ are contracted by the birational morphism $\alpha:X\to \bar{\mfk{M}}_{X,\mcl{F}}$.
\end{rem}

\vspace{3mm}

The rest of this subsection is devoted to the proof of the following theorem.

\begin{thm}\label{main'}
The morphism $\alpha:X\to \bar{\mfk{M}}_{X,\mcl{F}}$ is an isomorphism to the normalization of a fine moduli space of $G$-constellations whose universal family is pulled back by $\alpha$ to $\mcl{F}$ if and only if all the $G$-constellations in $\mcl{F}$ are indecomposable.
\end{thm}

We start with the following lemma.

\begin{lem}\label{lem:simplicial}
Suppose that a gnat family $\mcl{F}$ on a crepant resolution $X$ of $\C^n/G$ satisfies the equivalent conditions (a),(b),(c),(d) in Lemma \ref{closed}. Then $\bar{\mfk{M}}_{X,\mcl{F}}$ is $\Q$-factorial.
\qed
\end{lem}

{\em Proof.} For any exceptional divisor $E_\ell\subset X=X(\Sigma)$, the fan of $\bar{\mcl{V}}_{\Sigma,\mcl{F}}$ contains a ray $\tau_{\mcl{F},\ell}$ which is mapped under $q_G$ to the ray $\tau_\ell\in\Sigma$ for $E_\ell$ by Lemma \ref{crepant divisor}. Since $\mcl{F}$ satisfies (b), the morphism $q_2:\bar{\mcl{V}}_{\Sigma,\mcl{F}}\to \bar{\mfk{M}}_{X,\mcl{F}}$ is a geometric quotient. Thus, the set of the rays of the fan of $\bar{\mfk{M}}_{X,\mcl{F}}$ is the same as that of $\Sigma$ and has cardinality $n+m$.

We take these $n+m$ rays and any $(r-1)-m$ rays of $\Sigma_\mcl{V}$ (if $r-1>m$) so that these rays span $(N_\mcl{V})_\R$ as a vector space. Let $P=\C[y_1,\dots,y_{n+r-1}]$ be the polynomial ring whose variables are indexed by these rays. Then the lattice $\tilde{N}:=\Z^{n+r-1}$ corresponding to the toric variety $\mathrm{Spec}\,P$ admits a homomorphism to $N_\mcl{V}$ which sends the basis to the generators of the $n+r-1$ rays. By composing this homomorphism with $q_G$, we obtain a map $\tilde{N}\to N_G$ and we denote its kernel by $K$. Let $T_K\subset (\C^*)^{n+r-1}$ be the subtorus associated to the inclusion $K\subset\tilde{N}$. Then, by Proposition \ref{quot}, $\bar{\mfk{M}}_{X,\mcl{F}}$ is obtained as a good quotient of a certain open subset of $\mathrm{Spec}\,P$ by the $T_K$-action. It follows from the equality $\dim T_K=\dim \bar{T}^{r-1}$ that this quotient is in fact geometric. Since $\mathrm{Spec}\,P$ is smooth, the geometric quotient $\bar{\mfk{M}}_{X,\mcl{F}}$ of its open subset is $\Q$-factorial by Proposition \ref{quot}.
\qed

\vspace{3mm}

\begin{rem}
1. If $X\to \C^n/G$ is projective, then in fact $X$ is realized as the GIT-quotient of $\mathrm{Spec}\,P$ for a generic character $\chi$ of $T_K$. Although we obtain a generic stability condition $\theta$ via $\chi(T_K)_\R\cong\chi(\bar{T}^{r-1})_\R\cong \Theta$, this does not imply that $X\cong\bar{\mfk{M}}_\theta$. This is because a $\chi$-invariant function $f$ (namely, $f$ such that $t\cdot f=\chi(t)f$ for $t\in \bar{T}^{r-1}$) on $\mathrm{Spec}\,P$, which can be regarded as a $\chi$-invariant function on $\bar{\mcl{V}}_{\Sigma,\mcl{F}}$ after taking its powers,  might not extend to the whole $\bar{\mcl{V}}$.\\
2. In general $\bar{\mfk{M}}_\theta$ can be non-$\Q$-factorial even if we take a generic $\theta$. The above proof shows that $\bar{\mfk{M}}_\theta$ for generic $\theta$ is $\Q$-factorial if the number of irreducible exceptional divisors of $\bar{\mfk{M}}_\theta$ is less than $r$.
\qed
\end{rem}

\vspace{3mm}

Before we proceed, we give a general observation. Let $\sigma\in\Sigma$ be an $\ell$-dimensional cone. Then it is generated by the one-dimensional cones corresponding to some exceptional divisors $E_{k_1},\dots,E_{k_{\ell'}}$ and coordinate divisors $D_{k'_1},\dots,D_{k'_{\ell-\ell'}}$(cf. (\ref{delta})). By the construction of $\bar{\vphi}_\mcl{F}$, the associated cone $\tilde{\sigma}\in\tilde{\Sigma}$ is the smallest cone of the fan $\Sigma_{\mcl{V}}$ containing $\tau_{\mcl{F},{k_1}},\dots,\tau_{\mcl{F},k_{\ell'}}$ and $\tilde{e}_{k'_1},\dots,\tilde{e}_{k'_{\ell-\ell'}}$ where $\tilde{e}_j$ is the unique ray of $\Sigma_{\mcl{V}}$ satisfying $q_G(\tilde{e}_j)=\R_{\ge0} e_j$.

\vspace{3mm}

{\it Proof of Theorem \ref{main'}}\\
The ``only if" part follows since all $\bar{T}^{r-1}$-orbits in $\mcl{V}_{\Sigma,\mcl{F}}$ must be closed in order for $\bar{\mfk{M}}_{X,\mcl{F}}$ to be a fine moduli space (and hence to be a geometric quotient of $\bar{\mcl{V}}_{\Sigma,\mcl{F}}$).

We next consider the ``if" part. Assume that the four equivalent conditions (a), (b), (c), (d) in Lemma \ref{closed} are satisfied. For each $n$-dimensional cone $\sigma\in\Sigma$, we see that $q_G(\tilde{\sigma})$ contains $\sigma$ by the observation above, and thus $\tilde{\sigma}$ has dimension at least $n$. By toric geometry, the dimension of $\tilde{\sigma}$ is equal to the dimension of the stabilizer group $\mathrm{Stab}_{x_{\tilde{\sigma}}}\subset T_\mcl{V}$ where $x_{\tilde{\sigma}}$ is the distinguished point of the affine toric variety $X(\tilde{\sigma})\subset \bar{\mcl{V}}_{\Sigma,\mcl{F}}$. By the condition (c), the intersection of this stabilizer group and the subtorus $\bar{T}^{r-1}\subset T_\mcl{V}$ has dimension zero, which implies that $\dim(\tilde{\sigma})=n$. Therefore, $\tilde{\sigma}$ is an $n$-dimensional simplicial cone by Lemma \ref{lem:simplicial}. 

By the construction of the fan of $\bar{\mfk{M}}_{X,\mcl{F}}$, the projection $q_G$ maps $\tilde{\sigma}$ isomorphically to $\sigma$, and all of the four toric varieties in the diagram (\ref{morph}) has isomorphic associated fans. This particularly shows that $\alpha:X\to \bar{\mfk{M}}_{X,\mcl{F}}$ is an isomorphism. Since $\bar{\mfk{M}}_{X,\mcl{F}}$ is a geometric quotient by the condition (a), it is a fine moduli space by the argument of the proof of \cite[Proposition 5.3]{Kin} and its universal family is identified with $\mcl{F}$ via the isomorphism $\alpha$.
\qed

\vspace{3mm}

\begin{rem}
The results in this section still hold even if we replace crepant resolutions by relative minimal models. However, one can in fact show that every relative minimal model satisfying the four equivalent conditions in Lemma \ref{closed} is automatically smooth, and thus we only consider crepant resolutions.
\qed
\end{rem}

\vspace{3mm}

As mentioned in the introduction, giving a description of a crepant resolution as a fine moduli space is deeply related to the derived McKay correspondence. In dimension 2 or 3, the moduli space $\bar{\mfk{M}}_\theta$ with generic $\theta$ is always a crepant resolution of $\C^n/G$, and the universal family on $\bar{\mfk{M}}_\theta$ induces an equivalence between the triangulated categories $D^b(\bar{\mfk{M}}_\theta)$ and $D^b_G(\C^n)$ (\cite{BKR}) where $D^b(\bar{\mfk{M}}_\theta)$ and $D^b_G(\C^n)$ are the bounded derived category of coherent sheaves on $\bar{\mfk{M}}_\theta$ and that of $G$-equivariant sheave on $\C^n$ respectively. Then it is natural  to ask the following:

\begin{que}\label{question}
Let $G\subset SL_n(\C)$ be a finite abelian subgroup and let $X\to \C^n/G$ be a not-necessarily-projective crepant resolution. Suppose that a gnat family $\mcl{F}$ on $X$ parametrizes indecomposable $G$-constellations (and thus $X$ is regarded as a fine moduli space of $G$-constellations by Theorem \ref{main}). Then is it true that the Fourier-Mukai transform $\Phi_\mcl{F}:D^b(X)\to D^b_G(\C^n)$ by $\mcl{F}$ always gives an equivalence of derived categories?
\end{que}

The author does not know the answer to this question even in the case $n=3$. In \cite{L2}, a condition for $\Phi_\mcl{F}$ to be an equivalence is given in terms of fibers of $\mcl{F}$. However, the author could not find a clear relationship between this condition and the condition for $X$ to be a fine moduli space (that is, the indecomposability of the fibers of $\mcl{F}$).

\section{Special relative minimal models as coarse moduli spaces of $G$-constellations}\label{6}

Before stating the claim, let us recall that every pair of a gnat family $\mcl{F}$ on a relative minimal model $X$ and an exceptional divisor $E_k\subset X$ defines a torus-invariant divisor of $\mcl{V}$ and hence its corresponding one-dimensional cone $\tau_{\mcl{F},E_k}$ of the fan of $\bar{\mcl{V}}$ (cf. Lemma \ref{crepant divisor}). We will prove the following.

\begin{prop}\label{prop:special}
Let $G\subset SL_n(\C)$ be a finite abelian subgroup and let $X=X(\Sigma)$ be a relative minimal model obtained by repeated blowing-ups of $\C^n/G$. Then there is a gnat family $\mcl{F}$ on $X$ such that the projection $q_G:(N_\mcl{V})_\R\to (N)_\R$ in (\ref{diag1PS}) induces a good quotient $q_2:\bar{\mcl{V}}_{\Sigma,\mcl{F}}\to X$ (i.e $q_G^{-1}(\sigma)\cap |\tilde{\Sigma}|\in \tilde{\Sigma}$ for any $\sigma\in \Sigma$) and also satisfies $q_G^{-1}(\tau_k)\cap |\tilde{\Sigma}|=\tau_{\mcl{F},E_k}$ for all $k$.
\end{prop}

\begin{rem}
By using a similar argument in \cite[Proposition 5.2]{Kin}, we see that $X$ is regarded as a coarse moduli space of the $G$-constellations parametrized by $\mcl{F}$. 
\qed
\end{rem}

\vspace{3mm}

In terms fans of toric varieties, $X$ is obtained by applying star subdivisions of a simplicial fan $\Sigma'$ in $(N_G)_\R=\R^n$ repeatedly at the lattice points on the simplex $\{(u_i)\in\R_{\ge0}^n\mid \sum u_i=1\}$. Precisely, the {\it star subdivision} of $\Sigma'$ at a lattice point $v\in N_G\cap |\Sigma|$ is the minimal refinement $\Sigma_v$ of $\Sigma$ containing the simplicial cones
$$\sigma_\ell:=\mrm{Cone}(u_1,\dots,\hat{u}_\ell,\dots,u_n,v)$$
for all $n$-dimensional $\sigma=\mrm{Cone}(u_1,\dots,u_n)\in\Sigma'$ containing $v$ and all $\ell=1,\dots,n$. For example, if $v$ is inside the interior of an $n$-dimensional cone $\sigma\in \Sigma'$, then $\Sigma_v$ is obtained by removing $\sigma$ and adding $n$ simplicial cones $\sigma_\ell$ and their faces to $\Sigma'$. If we start with the fan $\Sigma_0$ containing $\sigma_0^G$ as the unique maximal cone and apply star subdivisions repeatedly at all the points of $\{v_g \mid g\in G:\text{junior}\}$, we obtain a relative minimal  model regardless of the order of the choices of $v_g$ (see Subsection \ref{3.2}). We fix the order $g_1,\dots,g_m$ of the junior elements of $G$ (or equivalently, the order of the exceptional divisors $E_k$ of $X$) and let $X_k$ be the partial resolution of $\C^n/G$ obtained by applying the star subdivisions, starting from $\Sigma_0$, repeatedly at the points $v_{g_1},\dots,v_{g_k}$ in this order. We choose the order of $g_k$'s so that $X=X_m$. We simply write the fan of $X_k$ as $\Sigma_k$.

\vspace{3mm}

{\it Proof of Proposition \ref{prop:special}}\\
We will choose a gnat family $\mcl{F}$ on $X=X_m$ with the desired property by induction on the order of $G$. We first choose the ``canonical family" $\mcl{F}_1$ on $X_1$ in the sense of \cite[\S3.3]{L1}. This is characterized as a unique gnat family whose coefficients $b_{i,1}$ satisfy $0\le b_{i,1}<1$. We assume that the generator $v$ of $\tau_1$ lies in the interior of $|\Sigma_0|$ for simplicity since the other cases are similar. Then the fan of $X_1$ is the star subdivision $(\sigma_0)_v$, which contains exactly $n$ maximal cones $\sigma_\ell:=(\sigma_0^G)_\ell,\,\ell=1,\dots,n$ defined as above. The distinguished $G$-constellation  $(\mcl{F}_1)_{\sigma_\ell}=\{p_{i,j}\}$ corresponding to $\sigma_\ell$ is obtained by
\begin{equation}\label{canonical}
p_{i,j}=
\begin{cases}
1&\text{ if }j=\ell\text{ and }b_{i,\ell}+v_\ell=b_{{k(i,\ell)},\ell}\\
0&\text{ otherwise}
\end{cases}
\end{equation}
where $v=\sum_j v_j e_j\,(v_j\in \Q)$.

Let $\tilde{\sigma}_\ell$ be the cone of the fan $\tilde{\Sigma}$ in $(N_\mcl{V})_\R$ corresponding to $(\mcl{F}_1)_{\sigma_\ell}$ or, in other words, the smallest cone containing $\tau_{\mcl{F}_1,E_k}$ and $\tilde{e}_1,\dots,\tilde{e}_{\ell-1},\tilde{e}_{\ell+1},\dots,\tilde{e}_n$ where $\tilde{e}_j$ is the ray of $\tilde{\Sigma}$ corresponding to the coordinate divisor $D_j\subset X$ (i.e. $q_G(\tilde{e}_j)=\R_{\ge0} e_j$). As we will see below, in fact $q_G(\tilde{\sigma}_\ell)=\sigma_\ell$ holds.

For each $\ell$, let $N_\ell\subset N_G$ be the sublattice generated by $v$ and $e_1,\dots,\hat{e}_\ell,\dots,e_n$. We then define a finite abelian group $G_\ell:=N_G/N_\ell$ so that we have an identification $N_{G_\ell}\to N_G$ that sends $e_i$ to $e_i$ for $i\ne\ell$ and $e_\ell$ to $v$ (cf. \cite[\S3]{J}). The basis $\xi_1,\dots,\xi_n$ of $M_\ell:=N_\ell^\vee\subset M\otimes_\Z \Q$ dual to $e_1,\dots,e_{\ell-1},v,e_{\ell+1},e_n$ is written as
$$\xi_j=
\begin{cases}
x_j x_\ell^{-\frac{v_j}{v_\ell}}&\text{ if }j\ne\ell\\
x_\ell^{\frac{1}{v_\ell}}&\text{ if }j=\ell
\end{cases}$$
where the monomials $x_1,\dots, x_n$ are regarded as the dual basis of $e_1,\dots,e_n\in N$. Then $G_\ell$ naturally acts on $\C^n=\mrm{Spec}\,\C[\xi_1,\dots,\xi_n]$ and thereby we regard $G_\ell$ as a subgroup of $SL_n(\C)$. Note that we have $M_G=(M_\ell)_{G_\ell}$. Following \cite{DLR} and \cite{J}, we define a (non-additive) function $\phi_\ell:M\to M_\ell$ (called the {\it round down function}) by
$$\phi_\ell(x_1^{m_1}\cdots x_n^{m_n})=\xi_1^{m_1}\cdots \xi_\ell^{\lfloor \sum_j v_j m_j \rfloor}\cdots \xi_n^{m_n}.$$
This also gives a surjection $\chi(G)\to \chi(G_\ell)$, which is also denoted by $\phi_\ell$ by abuse of notation, such that $\phi_\ell(\rho)\in \chi(G_\ell)$ is the character defined by the monomial $\phi_\ell(f)$ for any monomial $f\in M$ which defines the character $\rho\in \chi(G)$. We denote by $\chi'_j\in \chi(G_\ell)$ the character of $\xi_j$. 

Using $\phi_\ell$, we can lift any distinguished $G_\ell$-constellation to a distinguished $G$-constellation. To explain this, we should first note that $\phi_\ell(\rho_{k(i,j)})(=\phi_\ell(\rho_i\otimes \chi_j))$ may not coincide with $\phi_\ell(\rho_i)\otimes \chi'_j$ since $\phi_\ell$ is not multiplicative. From the definition of $\phi_\ell$, we see that we have the following three cases:\\
(1) $\phi_\ell(\rho_{k(i,j)})=\phi_\ell(\rho_i)\otimes \chi'_j$\\
(2) $\phi_\ell(\rho_{k(i,j)})=\phi_\ell(\rho_i)\otimes \chi'_\ell \chi'_j$\\
(3) $j=\ell$ and $\phi_\ell(\rho_{k(i,j)})=\phi_\ell(\rho_i)$.\\
Set $\mrm{Irr}(G_\ell)=\{\rho'_0,\dots,\rho'_{r'-1}\}\,(r':=|G_\ell|)$ with fixed basis of $\rho'_i$'s, and suppose we are given a distinguished $G_\ell$-constellation $F'$, namely, the regular representation $\bigoplus_{\rho'\in \mrm{Irr}(G_\ell)} \rho'$ admitting $G_\ell$-equivariant actions of $\xi_1,\dots,\xi_n$ such that the action map $\rho'\to\rho'\otimes \chi'_j$ by $\xi_j$ sends the fixed basis of $\rho'$ to that of $\rho'\otimes \chi_j$ or to zero. We define a $G$-constellation $\phi_\ell^*(F')=\{p_{i,j}\}_{i,j}$ by setting $p_{i,j}=0$ if\\
the case (1) holds and the action map $\phi_\ell(\rho_i)\to \phi_\ell(\rho_i)\otimes \chi'_j$ by $\xi_j$ is zero, or\\
the case (2) holds and the action map $\phi_\ell(\rho_i)\to \phi_\ell(\rho_i)\otimes \chi'_\ell\chi'_j$ by $\xi_\ell\xi_j$ is zero,\\
and otherwise $p_{i,j}=1$. One can see that $\phi_\ell^*(F')$ is actually a $G$-constellation by checking the relations $p_{k(i,j),j'}p_{i,j}-p_{k(i,j'),j}p_{i,j'}$.

\begin{rem}
The natural inverse $\phi_\ell^*$ defined in \cite[\S3.1]{J} is regarded as the special case of our $\phi_\ell^*$ above. Our construction does not require that the $G_\ell$-constellation $F'$ comes from a $G_\ell$-brick. 
\qed
\end{rem}

\begin{lem}\label{lifting}
Let $\mfk{F}_\ell$ be the set of distinguished $G$-constellations corresponding to faces of the cone $\tilde{\sigma}_\ell$. Then the assignment $F'\mapsto \phi_\ell^*(F')$ gives a bijection between the set of distinguished $G_\ell$-constellations and $\mfk{F}_\ell$.
\end{lem}

{\em Proof.} We construct an inverse $(\phi_\ell)_*$ of $\phi_\ell^*$. To do this, note that, for any $\rho'\in \chi(G_\ell)$ and $j$, there is $\rho_i\in \chi(G)$ such that $\phi_\ell(\rho_i)=\rho'$ and $\phi_\ell(\rho_{k(i,j)})=\phi_\ell(\rho_i)\otimes \chi'_j$ (cf. \cite[Lemma 3.9]{J}). Then, for each $F\in\mfk{F}_\ell$, we can define a distinguished $G_\ell$-constellation $(\phi_\ell)_*(F)$ such that its action map $\phi_\ell(\rho_i)\to\phi_\ell(\rho_i)\otimes \chi'_j$ by $\xi_j$ is nonzero if and only if the action map $\rho_i\to \rho_{k(i,j)}$ in $F$ by $x_j$ is nonzero. We show $\phi_\ell^*\circ(\phi_\ell)_*=id$. Showing $(\phi_\ell)_*\circ\phi_\ell^*=id$ can be done in a similar spirit. 

We take any $F\in\mfk{F}_\ell$. In order to show $\phi_\ell^*((\phi_\ell)_*(F))=F$, it suffices to check that the actions $\rho_i\to \rho_{k(i,j)}$ by $x_j$ in the both sides coincide for each $i,j$. We first consider the case when $p_{i,j}=1$ for $F=\{p_{i,j}\}_{i,j}$. If (1) $\phi_\ell(\rho_{k(i,j)})=\phi_\ell(\rho_i)\otimes \chi'_j$ holds, then the action $\phi_\ell(\rho_i)\to\phi_\ell(\rho_i)\otimes \chi'_j$ by $\xi_j$ for $(\phi_\ell)_*(F)$ is nonzero and thus we have $p'_{i,j}=1$ for $\phi_\ell^*((\phi_\ell)_*(F))=\{p'_{i,j}\}_{i,j}$. If (2) $\phi_\ell(\rho_{k(i,j)})=\phi_\ell(\rho_i)\otimes \chi'_\ell\chi'_j$ holds, then one can show that both of the actions
$$\phi_\ell(\rho_i)\to\phi_\ell(\rho_i)\otimes \chi'_\ell$$
by $\xi_\ell$ and
$$\phi_\ell(\rho_{k(i,\ell)})\to\phi_\ell(\rho_{k(i,\ell)})\otimes \chi'_j$$
by $\xi_j$ for $(\phi_\ell)_*(F)$ are nonzero and hence $p'_{i,j}=1$. This follows from $p_{k(i,j),\ell}=1$ and in particular $p_{i,\ell}=p_{k(i,\ell),j}=1$. Note that the the assumption $F\in\mfk{F}_\ell$ implies that the action map $\rho_i\to \rho_{k(i,\ell)}$ in $F$ by $x_\ell$ is nonzero as long as the action map $\rho_i\to \rho_{k(i,\ell)}$ in $(\mcl{F}_1)_{\sigma_\ell}$ by $x_\ell$ is nonzero (cf. (\ref{canonical})). Finally if (3) $j=\ell$ and $\phi_\ell(\rho_{k(i,j)})=\phi_\ell(\rho_i)$ hold, then we clearly have $p'_{i,j}=1$.

Next we consider the case when $p_{i,j}=0$. In the case (1), $p_{i,j}=0$ follows immediately from the definitions of $(\phi_\ell)_*$ and $\phi_\ell^*$. In the case (2), one can show $p'_{i,j}=0$ from $p_{i,\ell} p_{k(i,\ell),j}=p_{i,j}p_{k(i,j),\ell}=0$. The case (3) does not happen since $F\in \mfk{F}_\ell$. Therefore, in any case we have $p_{i,j}=p'_{i,j}$ and the claim follows.
\qed

\vspace{3mm}

Let $\mcl{V}_\ell$ be the space of $G_\ell$-constellations (i.e. the counterpart of $\mcl{V}$ for $G_\ell$ instead of $G$). Similarly, we can consider the torus $\bar{T}^{r'-1}$ for $G_\ell$, and we have a natural embedding $\iota:\bar{T}^{r'-1}\to \bar{T}^{r-1}$ such that $\iota(\bar{t_{i'}})$ acts on $\rho_i$ as multiplication of $t_i$ (resp. trivially) if $\phi_\ell(\rho_i)=\rho_{i'}$ (resp. $\ne\rho_{i'}$). This also induces a homomorphism $T_{\mcl{V}_\ell}\to T_\mcl{V}$ of the big tori. By the correspondence of $G_\ell$- and $G$-constellations in Lemma \ref{lifting}, we obtain a  homomorphism $(N_{\mcl{V}_\ell})_\R\to (N_\mcl{V})_\R$ which sends the fan $\Sigma_\ell$ of $\mcl{V}_\ell$ to the cone $\tilde{\sigma}_\ell$ (regarded as a fan) isomorphically. Moreover, $N_{\mcl{V}_\ell}\to N_\mcl{V}$ is compatible with the projections $q_{G_\ell}:(N_{\mcl{V}_\ell})_\R\to (N_{G_\ell})_\R\cong(N_G)_\R$ and $q_G:(N_{\mcl{V}})_\R\to (N_G)_\R$, which follows from the fact that $\phi_\ell^*(t\cdot F')=\iota(t)\cdot \phi_\ell^*(F')$ for any $G_\ell$-constellation $F'$ and $t\in \bar{T}^{r'-1}$.

By the assumption of $X=X_m$, it admits an open subset $Y_\ell\subset X$ which is obtained by applying star subdivisions to $\sigma_\ell$. Note that $Y_\ell$ is also a relative minimal model of $\C^n/G_\ell$. By induction on $|G|$, we have a gnat family $\mcl{F}_{Y_\ell}$ of $G_\ell$-constellations on $Y_\ell$ satisfying the conditions in Proposition \ref{prop:special} for $G_\ell$ instead of $G$. The gnat families $\mcl{F}_{Y_\ell}$ on $Y_\ell$ for $G_\ell$ are lifted by $\phi_\ell^*$ to be gnat families on $Y_\ell\subset X$ for $G$ by Lemma \ref{lifting} and then glue to become a gnat family $\mcl{F}$ on $X$. We immediately see that $\mcl{F}$ satisfies the desired properties since each $\mcl{F}_{Y_\ell}$ satisfies the same properties. This completes the proof of Proposition \ref{prop:special}.\qed

\vspace{3mm}

\begin{rem}
Since $X=X_m$ is projective over $\C^n/G$, one may ask if $X$ is realized as $\bar{\mfk{M}}_\theta$ for some (possibly non-generic) $\theta\in \Theta$. The author does not know the answer but this seems true at least for small values of $n$ and $r$. By induction, it suffices to show $X_1\cong \bar{\mfk{M}}_{\theta'}$ for some $\theta'$ in order to show $X\cong \bar{\mfk{M}}_{\theta}$. As the referee suggests, one might be able to answer the question by making good use of \cite[Algorithm 7.6]{CMT}. This algorithm gives an explicit procedure to determine the set of stability conditions $\theta$ for which a given $G$-constellation is $\theta$-(semi)stable.
\qed
\end{rem}

\vspace{3mm}

\begin{ex}\label{A4}($A_4$-singularity)
We use similar notation to Example \ref{A3}. Here, we let $G\subset SL_2(\C)$ be the cyclic group of order 5 so that $\C^2/G$ has $A_4$-singularity. Let $X_1$ be the toric variety obtained by the star subdivision of $\sigma_0=\mrm{Cone}(e_1,e_2)$ at $v=\frac{3}{5}e_1+\frac{2}{5}e_2$, and let $E_1\subset X_1$ be the exceptional divisor. The canonical family $\mcl{F}_1$ on $X_1$ is obtained by setting
$$D_{\rho_0}=0,\;D_{\rho_1}=\frac{3}{5}E_1,\;D_{\rho_2}=\frac{1}{5}E_1,\;D_{\rho_3}=\frac{4}{5}E_1,\;D_{\rho_4}=\frac{2}{5}E_1.$$
Then the $G$-constellation $\mcl{F}_1(\sigma_1)$ with $\sigma_1:=(\sigma_0^G)_1=\mrm{Cone}(v,e_2)$ admits three indecomposable summands, which are equal to $\rho_0\oplus \rho_1,\,\rho_2\oplus \rho_3$ and $\rho_4$ as $G$-modules respectively. Note that $y$ acts trivially on $\mcl{F}_1(\sigma_1)$.

The group $G_1$ for $\sigma_1$ is a cyclic group of order 3 and thus $\sigma_1$ is identified with the cone $\sigma_0^{G_1}$ for $A_2$-singularity. We take the canonical family on the minimal resolution of $\C^2/G_1$ and let $F'$ be its distinguished $G_1$-constellation for the exceptional divisor $E'$ corresponding to the ray $\R_{\ge0}(\frac{2}{3}v+\frac{1}{3}e_2)\subset \sigma_0^{G_1}$. One can compute $\phi_1^*(F')$ by noticing
$$\phi_1^{-1}(\rho'_0)=\{\rho_0,\rho_1\},\phi_1^{-1}(\rho'_1)=\{\rho_2,\rho_3\},\phi_1^{-1}(\rho'_2)=\{\rho_4\}$$
where $\rho'_i\in \chi(G_1)$ is the $i$-th power of the character determined by $\xi_1=x^{\frac{5}{3}}$. We see that $\phi_1^*(F')$ is obtained by adding the nontrivial actions $\rho_1\to \rho_2$ by $x$ and $\rho_0\to \rho_4$ and $\rho_4\to \rho_3$ by $y$ to $\mcl{F}_1(\sigma_1)$.

The lifted $G$-constellation $\phi_1^*(F')$ is a distinguished one for the exceptional divisor $E_2\subset X$ corresponding to the ray $\R_{\ge0}(\frac{2}{3}v+\frac{1}{3}e_2)=\R_{\ge0}(\frac{2}{5}e_1+\frac{3}{5}e_2)\subset \sigma_0^G$ where $X$ is the minimal resolution of $\C^2/G$. However, one can check that $\phi_1^*(F')$ is not the $G$-constellation for the canonical family (but for the ``maximal shift family" in the sense of \cite[\S3.5]{L1}). Therefore, the resulting gnat family $\mcl{F}$ on the resolution $X$ constructed as above is different from the canonical family as a whole in general.
\qed
\end{ex}

\begin{center}
National Center for Theoretical Sciences, Taipei, Taiwan\\
Email address: \texttt{yamagishi@ncts.ntu.edu.tw}
\end{center}

\end{document}